\documentclass[a4paper]{amsart}
\usepackage{amsmath,amssymb,amsfonts}
\usepackage{mathrsfs,latexsym,amsthm,stmaryrd,enumerate,epic}
\usepackage{amscd}  
\usepackage[all]{xy}

\date{\today}

\newtheorem{theorem}{Theorem}
\newtheorem{lemma}[theorem]{Lemma}
\newtheorem{corollary}[theorem]{Corollary}
\newtheorem{proposition}[theorem]{Proposition}

\newtheorem{condition}[theorem]{Condition}
\newtheorem{remark}[theorem]{Remark}

\renewcommand{\phi}{\varphi}

\newcommand{\ob}{\mathrm{Ob}\,}
\newcommand{\mi}{\mathtt{i}}
\newcommand{\mj}{\mathtt{j}}
\newcommand{\mk}{\mathtt{k}}
\newcommand{\ml}{\mathtt{l}}
\newcommand{\mm}{\mathtt{m}}
\newcommand{\mn}{\mathtt{n}}

\newcommand{\mt}{\mathtt{t}}
\newcommand{\mone}{\mathtt{1}}

\title[Partialization of categories]{Partialization of 
categories and inverse braid-permutation monoid}

\author{Ganna Kudryavtseva}
\address{G. K.: Department of Mechanics and Mathematics, Kyiv
Taras Shevchenko University, 64 Volodymyrska st., UA-01033, Kyiv,
UKRAINE.}
\email{akudr\symbol{64}univ.kiev.ua}
\thanks{G. K. was supported by the Royal Swedish Academy of 
Sciences}
\author{Volodymyr Mazorchuk}
\address{V. M.: Department of Mathematics, Uppsala University, 
Box 480, SE-75106, Uppsala, SWEDEN.}
\email{mazor\symbol{64}math.uu.se}
\thanks{V. M. was supported by STINT, the Royal Swedish Academy of 
Sciences, the Swedish Research Council.}
\numberwithin{equation}{section}

\begin{document}

\begin{abstract}
We show how the categorial approach to inverse monoids 
can be described as a certain endofunctor (which we 
call the partialization functor) of some category. 
In this paper we show that this functor can be used to 
obtain several recently defined inverse monoids, and
use it to define a new object, which we call
the inverse braid-permutation monoid. A presentation 
for this monoid is obtained. Finally, we study some 
abstract properties of the partialization functor and
its iterations. This leads to a categorification of a
monoid of all order preserving maps, and series of 
orthodox generalizations of the symmetric inverse semigroup.
\end{abstract}

\maketitle

\tableofcontents

\section{Introduction}\label{s1}

The study of inverse semigroups forms a classical part of
the theory of semigroups. Several approaches to construction
of inverse semigroups are known. The most abstract one seems
to be  the categorical approach, which was worked out in various
variations in \cite{Ho,Le,La}. In Section~\ref{s2} we
give a short overview of this approach in the interpretation,
most suitable for our further purposes. Roughly speaking,
starting with a small category with pullbacks for monomorphisms,
there is a functorial way to enlarge the set of morphisms of
this category by what is naturally to call ``partial morphisms''.
This defines a functor, which we call the {\em partialization
functor}. Under some mild assumptions the endomorphism monoids of the
new category turn out to be inverse monoids. For example,
starting from the category of all finite sets in which the 
morphisms are all monomorphisms, the described above
partialization procedure gives us a new category, where the endomorphism
monoids are just symmetric inverse monoids. Several other examples
are discussed in Section~\ref{s3}.

One can apply this approach 
to produce several inverse monoids which recently appeared in 
the literature. The main example which we have in mind is the 
inverse braid monoid, defined and studied in \cite{EL}.
The idea is the following: the braid group can be realized as
the mapping class group (of homeomorphisms with compact support) 
of a punctured plane. In Section~\ref{s4} we define the category,
whose objects are punctured planes with different punctures
and whose morphism sets are isotopy classes of homeomorphisms
between such planes. It turns out that the procedure of 
categorical partialization is directly applicable. The result is
a new category, in which the endomorphism monoids are exactly
the inverse braid monoids from \cite{EL}.

The latter example has an immediate generalization. Instead of the
punctured plane let us consider the space of $n$ unknotted
and unlinked simple closed oriented curves in $\mathbb{R}^3$. 
The corresponding motion group is the so-called {\em braid-permutation
group}, defined in \cite{Da} and studied in \cite{Wa,Ru,Mc,FRR}.
The categorical partialization is again directly applicable
and produces a category, in which it is natural to call endomorphism 
monoids the {\em inverse braid-permutation monoids}. These monoids
are new. We describe all constructions necessary for its
definition in Section~\ref{s5}, where we also obtain a 
presentation for this monoid.

Finally, in Section~\ref{s6} we go back to the abstract study 
of the partialization functor. We show that this functor is in
fact an endofunctor of a certain category. We explicitly describe
its iterations and show how they are connected with each other
by some canonical natural transformations. In this way we obtain a
categorification of the monoid of all order preserving maps on a chain,
which fix the endpoints. We also describe certain quasi-iterations of 
the partialization functor on a somewhat bigger category.
In natural examples the iterations of the partialization  functor 
produce orthodox monoids. We pay special attention to orthodox 
generalizations of the symmetric inverse monoid.

For a monoid, $S$ we denote by $E(S)$ the set of all idempotents
of $S$ and by $G(S)$ the group of units of $S$. Green's relations
on $S$ are denoted by  $\mathcal{L}$, $\mathcal{R}$, 
$\mathcal{H}$, $\mathcal{D}$, and $\mathcal{J}$.
Recall that a monoid $S$ is called {\em factorizable} if 
$S=E(S)G(S)$. 

\subsection*{Acknowledgments}

We would like to thank Ryszard Rubinsztein for many helpful
discussions.

\section{Partializations of a category}\label{s2}

In this section we recall the partialization procedure for
small categories, worked out in \cite{Ho,Le,La}. However,
our setup is slightly different from the one available
in the literature. It is more adapted to our further needs.

\subsection{The partialization functor $\mathscr{P}$}\label{s2.1}

In what follows $\mathscr{C}$ is a small category. For 
$\mi,\mj\in\mathscr{C}$ we denote by $\mathscr{C}(\mi,\mj)$ the 
set of morphisms from  $\mi$ to $\mj$. For $\mi\in\mathscr{C}$ 
we denote by  $\mone_{\mi}$ the identity morphisms in 
$\mathscr{C}(\mi,\mi)$.

Assume that $\mathscr{C}$ satisfies the following condition:
\begin{condition}\label{pbm}
For any $\mi,\mj,\mk\in \mathscr{C}$, any
$f\in \mathscr{C}(\mi,\mj)$, and any monomorphism 
$\alpha\in \mathscr{C}(\mk,\mj)$ there exists a pullback diagram
\begin{displaymath}
\xymatrix{
\mi\ar[r]^{f}&\mj\\
\mm\ar@{.>}[r]^{\hat{f}}
\ar@{^{(}.>}[u]^{\beta}&\mk\ar@{^{(}->}[u]_{\alpha}.
}
\end{displaymath}
\end{condition}

We are now going to define a new category, $\mathscr{P}(\mathscr{C})$,
as follows: $\mathscr{P}(\mathscr{C})$ has the same objects as
$\mathscr{C}$. For $\mi,\mj\in \mathscr{P}(\mathscr{C})$ to 
define $\mathscr{P}(\mathscr{C})(\mi,\mj)$ we
consider all possible diagrams of the form
\begin{equation}\label{eqeq1}
\mathrm{D}(\mi,\mj,\mk,\alpha,f):\quad\quad\xymatrix{
\mi&\mk\ar@{_{(}->}[l]_{\alpha} \ar[r]^{f}&\mj
}\quad\quad\text{ or }\quad\quad
\xymatrix{
\mi\\
\mk\ar@{^{(}->}[u]^{\alpha} \ar[r]^{f}&\mj.
}
\end{equation}
We will say that the diagrams $\mathrm{D}(\mi,\mj,\mk,\alpha,f)$ 
and $\mathrm{D}(\mi',\mj',\mk',\alpha',f')$  are equivalent 
provided that $\mi'=\mi$, $\mj'=\mj$, and there is  an isomorphism, 
$\gamma$, which makes the following diagram commutative:
\begin{equation}\label{neqeq2}
\xymatrix{
\mi\ar@{=}[d]&\mk\ar@{_{(}->}[l]_{\alpha} \ar[r]^{f}
\ar@/^/@{.>}[d]|-{\gamma}
&\mj\ar@{=}[d]\\
\mi&\mk'\ar@{_{(}->}[l]_{\alpha'}\ar@/^/@{.>}[u]|-{\gamma^{-1}} 
\ar[r]^{f'}&\mj.
}
\end{equation}
Denote by  $\mathscr{P}(\mathscr{C})(\mi,\mj)$ the set of all
equivalence classes. We define the composition of 
$\mathrm{D}(\mi,\mj,\mk,\alpha,f)$ with 
$\mathrm{D}(\mj,\mm,\mn,\beta,g)$ as
$\mathrm{D}(\mi,\mm,\ml,\alpha\gamma,gh)$ via the following
pullback diagram, whose existence is guaranteed by 
Condition~\ref{pbm}:
\begin{equation}\label{eqeq2}
\xymatrix{
\mi&&\\
\mk\ar@{^{(}->}[u]^{\alpha}\ar[r]^{f}&\mj&\\
\ml\ar@{^{(}.>}[u]^{\gamma}\ar@{.>}[r]^{h}&
\mn\ar@{^{(}->}[u]^{\beta}\ar[r]^{g}&\mm
}.
\end{equation}

It is straightforward to verify that the above composition is
well-defined and associative and hence $\mathscr{P}(\mathscr{C})$
is a category. 

For $\mi\in\mathscr{C}$ set $\mathfrak{i}(\mi)=\mi\in
\mathscr{P}(\mathscr{C})$; and for $f\in \mathscr{C}(\mi,\mj)$ set
$\mathfrak{i}(f)=\mathrm{D}(\mi,\mj,\mi,\mone_{\mi},f)$. 
This obviously defines a faithfull functor, 
$\mathfrak{i}:\mathscr{C}\to \mathscr{P}(\mathscr{C})$, which we
will call the {\em canonical inclusion}.

Let now $\mathscr{C}$ and $\mathscr{C}'$ be two categories satisfying
Condition~\ref{pbm} and let $F:\mathscr{C}\to \mathscr{C}'$ be a
functor, which preserves monomorphisms. Then $F$ maps the diagram
$\mathrm{D}(\mi,\mj,\mk,\alpha,f)$ to 
$\mathrm{D}(F(\mi),F(\mj),F(\mk),F(\alpha),F(f))$ and the
functorial properties of $F$ ensure that the equivalent diagrams
end up in equivalent diagrams. Hence $F$ induces a functor from
$\mathscr{P}(\mathscr{C})$ to $\mathscr{P}(\mathscr{C}')$, which
we will denote by $\mathscr{P}(F)$. In particular, 
$\mathscr{P}$ becomes a functor from the category of all small
categories, satisfying Condition~\ref{pbm}, where morphisms are
all functors, preserving monomorphisms, to the category of
all small categories. We will call $\mathscr{P}$ the
{\em (first) partialization functor}. The functor $\mathfrak{i}$, 
when applied to all possible categories satisfying
Condition~\ref{pbm}, gives a natural  transformation from  the 
canonical inclusion functor to  $\mathscr{P}$. 

\subsection{$\mathscr{P}$ and inverse semigroups}\label{s2.2}

From now on we always assume that $\mathscr{C}$ is
a small category satisfying Condition~\ref{pbm}. The 
connection of $\mathscr{P}$ to inverse semigroups is given 
by the following statement:

\begin{proposition}\label{pmonoinv}
Suppose that all morphisms in $\mathscr{C}$ are monomorphisms. 
Then $\mathrm{End}_{\mathscr{P}(\mathscr{C})}\big(\mi\big)$ is an 
inverse semigroup for every object $\mi\in \ob\mathscr{C}$.
\end{proposition}

\begin{proof}
We start the proof with the following easy observation, which one proves
by a direct calculation:

\begin{lemma}\label{lem:regularity}
Let $\mi,\mj,\ml\in\mathscr{C}$. Let further
$\alpha\in \mathscr{C}(\ml,\mi)$ and $\beta\in \mathscr{C}(\ml,\mj)$
be monomorphisms. Then the elements $x=\mathrm{D}(\mi,\mj,\ml,\alpha,\beta)$
and $y=\mathrm{D}(\mj,\mi,\ml,\beta,\alpha)$ satisfy $xyx=x$ and
$yxy=y$.
\end{lemma}

It follows from Lemma~\ref{lem:regularity} that $\mathrm{End}_{\mathscr{P}(\mathscr{C})}\big(\mi\big)$ is a regular 
semigroup. Hence, to prove our proposition we have just to show that the 
idempotents of $\mathrm{End}_{\mathscr{P}(\mathscr{C})}\big(\mi\big)$ 
commute. We proceed by describing the idempotents:

\begin{lemma}\label{lem:idempotents}
The element $x=\mathrm{D}(\mi,\mi,\ml,\alpha,\beta)\in
\mathscr{P}(\mathscr{C})(\mi,\mi)$ is an idempotent if and only
if $\alpha=\beta$.
\end{lemma}

\begin{proof}
If $\alpha=\beta$, then $x^2=x$ is checked by a direct calculation. 
Consider the pullback diagram
\begin{equation}\label{neqeq3}
\xymatrix{
\ml\ar@{^{(}->}[r]^{\beta}&\mi \\
\ml'\ar@{^{(}->}[u]^{\alpha'} \ar@{^{(}->}[r]^{\beta'}&\ml\ar@{^{(}->}[u]^{\alpha}.
}
\end{equation}
Then $x^2=\mathrm{D}(\mi,\mi,\ml',\alpha\alpha',\beta\beta')$.
Using \eqref{neqeq2}, the equivalence of $x^2$ and $x$ implies
the existence of an isomorphism $\gamma\in \mathscr{C}(\ml',\ml)$
such that $\alpha\gamma=\alpha\alpha'$ and  $\beta\gamma=\beta\beta'$. 
Since $\alpha$ and $\beta$ are monomorphisms, the latter equalities 
imply $\gamma=\alpha'=\beta'$, in particular, both $\alpha'$ 
and $\beta'$ are isomorphisms. Now \eqref{neqeq3} implies
$\alpha=\beta$.
\end{proof}

Let now $x=\mathrm{D}(\mi,\mi,\ml,\alpha,\alpha)$ and 
$y=\mathrm{D}(\mi,\mi,\ml',\beta,\beta)$ and 
\begin{displaymath}
\xymatrix{
\ml'\ar@{^{(}->}[r]^{\beta}&\mi \\
\mt\ar@{^{(}->}[u]^{\delta} \ar@{^{(}->}[r]^{\gamma}&\ml\ar@{^{(}->}[u]^{\alpha},
}
\end{displaymath}
be the pullback. Then $\alpha\gamma=\beta\delta$ and a direct 
calculation imply $xy=yx$. This completes the proof.
\end{proof}

\subsection{The partialization functor $\mathscr{Q}$}\label{s2.3}

Let $\mathscr{C}$ be as in Subsection~\ref{s2.1}. 
For $\mi,\mj\in \mathscr{C}$ denote by 
$\mathscr{Q}(\mathscr{C})(\mi,\mj)$ the subset of 
$\mathscr{P}(\mathscr{C})(\mi,\mj)$ consisting of all
$\mathrm{D}(\mi,\mj,\mk,\alpha,f)$ for which there exists
$\overline{f}\in \mathscr{C}(\mi,\mj)$ making the following
diagram commutative:
\begin{equation}\label{eqeq5}
\xymatrix{
\mi\ar@{.>}[rd]^{\overline{f}}\\
\mk\ar@{^{(}->}[u]^{\alpha} \ar[r]_{f}&\mj.
}
\end{equation}
It is easy to see that $\mathscr{Q}(\mathscr{C})$ is  a
subcategory of $\mathscr{P}(\mathscr{C})$ and that 
the natural inclusion $\mathscr{C}\hookrightarrow 
\mathscr{P}(\mathscr{C})$ factors through 
$\mathscr{Q}(\mathscr{C})$. From this in the same way as in
Subsection~\ref{s2.1} one gets that $\mathscr{Q}$ defines
a functor, which we will call the {\em (second) partialization 
functor}. The natural transformation $\mathfrak{i}$ then
factors as $\mathfrak{i}=\mathfrak{i}''\mathfrak{i}'$,
where $\mathfrak{i}'$ is the natural transformation from
the natural inclusion to $\mathscr{Q}$, and $\mathfrak{i}''$
is the natural transformation (inclusion)
from $\mathscr{Q}$ to $\mathscr{P}$.

The ``philosophical'' difference between $\mathscr{P}$
and $\mathscr{Q}$ is how one understands the notion of 
a {\em partial map}: for $\mathscr{P}$ a partial map is a map, 
defined on a subobject; whereas for $\mathscr{Q}$ a partial map 
is a restriction of an ordinary  map to a subobject. Both 
notions have their advantages and disadvantages, which could be
seen for example after comparing Subsection~\ref{s2.2} and
Subsection~\ref{s2.4}.

The element in \eqref{eqeq5} will be denoted by
$\overline{\mathrm{D}}(\mi,\mj,\mk,\alpha,f,\overline{f})$.
Note that $f=\overline{f}\alpha$ by definition.

\subsection{$\mathscr{Q}$ and inverse semigroups}\label{s2.4}

\begin{proposition}\label{thth8}
Suppose that  $\mathscr{C}$ is as above and 
$\mi\in \ob\mathscr{C}$ is such that
$\mathrm{End}_{\mathscr{C}}\big(\mi\big)$ is a group. 
Then  the semigroup  
$\mathrm{End}_{\mathscr{Q}(\mathscr{C})}\big(\mi\big)$ is 
an inverse monoid.
\end{proposition}

\begin{proof}
We will show that 
$\mathrm{End}_{\mathscr{Q}(\mathscr{C})}\big(\mi\big)$ is a
regular monoid with commuting idempotents.

\begin{lemma}\label{lele6}
$\mathrm{End}_{\mathscr{Q}(\mathscr{C})}\big(\mi\big)$ is 
regular.
\end{lemma}

\begin{proof}
Let $x=\overline{\mathrm{D}}(\mi,\mi,\mk,\alpha,f,g)$, where
$f=g\alpha$, be an element of $\mathscr{Q}(\mathscr{C})(\mi,\mi)$. 
Since $\mathrm{End}_{\mathscr{C}}\big(\mi\big)$ is a group, we can consider
$y=\overline{\mathrm{D}}(\mi,\mi,\mi,\mone_{\mi},g^{-1},g^{-1})$. 
A direct calculation shows  that $xyx=x$ and hence
$x$ is regular. This completes the proof.
\end{proof}

\begin{lemma}\label{lele7}
\begin{displaymath}
E(\mathrm{End}_{\mathscr{Q}(\mathscr{C})}\big(\mi\big))=
\{\overline{\mathrm{D}}(\mi,\mi,\mk,\alpha,\alpha,\mone_{\mi})\}.
\end{displaymath}
\end{lemma}

\begin{proof}
That each $\overline{\mathrm{D}}(\mi,\mi,\mk,\alpha,\alpha,\mone_{\mi})$ 
is an idempotent is obtained by a direct calculation.

Now let $x=\overline{\mathrm{D}}(\mi,\mi,\mk,\alpha,f,g)$,
$g=f\alpha$ be an idempotent. Then 
$\mathfrak{i}''(x)=\mathrm{D}(\mi,\mi,\mk,\alpha,f)$ is an 
idempotent as well. Since $g$ is an isomorphism, 
$f$ is a monomorphism and hence Lemma~\ref{lem:idempotents}
implies that without loss of generality we can assume
$\alpha=f$. The statement follows.
\end{proof}

Now a direct calculation shows that for 
$x,y\in \{\overline{\mathrm{D}}(\mi,\mi,\mk,\alpha,\alpha,\mone_{\mi})\}$
we have $xy=yx$. This completes the proof of the proposition.
\end{proof}

\begin{corollary}\label{cornew41}
Assume all conditions of Proposition~\ref{thth8}.
\begin{enumerate}[(i)]
\item\label{cornew41.1}
The monoid
$\mathrm{End}_{\mathscr{Q}(\mathscr{C})}\big(\mi\big)$ is
factorizable. 
\item\label{cornew41.2}
$\mathrm{End}_{\mathscr{Q}(\mathscr{C})}\big(\mi\big)$ is
a maximal factorizable subset of 
$\mathrm{End}_{\mathscr{P}(\mathscr{C})}\big(\mi\big)$.
\item\label{cornew41.3}
$\mathrm{End}_{\mathscr{Q}(\mathscr{C})}\big(\mi\big)=
\mathrm{End}_{\mathscr{P}(\mathscr{C})}\big(\mi\big)$
if and only if 
$\mathrm{End}_{\mathscr{P}(\mathscr{C})}\big(\mi\big)$
is factorizable.
\end{enumerate}
\end{corollary}

\begin{proof}
\eqref{cornew41.1} follows from the definition of $\mathscr{Q}$.
From the proofs of Proposition~\ref{pmonoinv}
and Proposition~\ref{thth8} we have
\begin{displaymath}
E(\mathrm{End}_{\mathscr{Q}(\mathscr{C})}\big(\mi\big))=
E(\mathrm{End}_{\mathscr{P}(\mathscr{C})}\big(\mi\big))
\quad\text{ and }\quad
G(\mathrm{End}_{\mathscr{Q}(\mathscr{C})}\big(\mi\big))=
G(\mathrm{End}_{\mathscr{P}(\mathscr{C})}\big(\mi\big)).
\end{displaymath}
This and \eqref{cornew41.1} imply \eqref{cornew41.2}, and
\eqref{cornew41.3} follows from \eqref{cornew41.2}.
\end{proof}

\subsection{Dual constructions}\label{s2.6}

All the constructions above admit obvious dualization
(in other words, the dualization reduces the dual constructions 
to the two constructions described above).
The semigroup-theoretical effect of this is that one
obtains the opposite semigroup (category). This might
look trivial, but in fact it is not. Natural categories
might have rather non-symmetric structure. Thus the
category of all finite sets, where morphisms are
all injections of sets, is {\bf not} dual to the
category of all finite sets, where morphisms are
all surjections of sets. The latter category satisfies 
the condition of pushouts for epimorphisms, and hence its
opposite category satisfies Condition~\ref{pbm}. Applied to 
these two very different categories, our construction
produces different inverse monoids. This will be discussed
later on in examples.

\section{Examples from the theory of semigroups}\label{s3}

In this section we show that many classical examples of
inverse semigroups can be obtained using
$\mathscr{P}$ or $\mathscr{Q}$.

\subsection{Finite symmetric inverse semigroup}\label{s3.1}

Let $\mathscr{C}_{1}$ denote the category, whose objects are finite
sets and morphisms are injective maps. In particular, all morphisms 
in this category are monomorphisms. If $X,Y,Z\in \mathscr{C}_{1}$,
$f:X\to Y$, $g:Z\to Y$, we can define $U=\{x\in X\,:\, f(x)\in g(Z)\}$ 
and for $u\in U$ define $h(u)$ as the unique element of
$Z$ such that $f(u)=g(h(u))$. The map $h$ is obviously an injection.
One shows that the commutative diagram:
\begin{displaymath}
\xymatrix{
X\ar@{^{(}->}[r]^{f} & Y \\
U\ar@{^{(}->}[u]^{incl}\ar@{^{(}->}[r]^{h} 
& Z\ar@{^{(}->}[u]^{g}
}
\end{displaymath}
is a pullback. Hence $\mathscr{C}_{1}$ satisfies
Condition~\ref{pbm}.  For $X\in \mathscr{C}_{1}$ we have that
$\mathrm{End}_{\mathscr{C}_{1}}\big(X\big)$ is the symmetric group on
$X$. A direct computation shows that 
$\mathrm{End}_{\mathscr{P}(\mathscr{C}_{1})}\big(X\big)=
\mathrm{End}_{\mathscr{Q}(\mathscr{C}_{1})}\big(X\big)$ is the full 
symmetric inverse semigroup 
$\mathcal{IS}(X)$ on $X$. Note that, although
$\mathrm{End}_{\mathscr{P}(\mathscr{C}_{1})}\big(X\big)=
\mathrm{End}_{\mathscr{Q}(\mathscr{C}_{1})}\big(X\big)$ for each
$X\in \mathscr{C}_{1}$, the categories
$\mathscr{P}(\mathscr{C}_{1})$ and 
$\mathscr{Q}(\mathscr{C}_{1})$ are different:
if $X,Y\in \mathscr{C}_{1}$ and $0<|X|<|Y|$, then 
$\mathscr{P}(\mathscr{C}_{1})(Y,X)\neq \varnothing$ while
$\mathscr{Q}(\mathscr{C}_{1})(Y,X)=\varnothing$.

\subsection{Dual symmetric inverse semigroup}\label{s3.15}

Let $\mathscr{C}_{2}$ denote the category, whose objects are finite
sets and morphisms are surjective maps. In particular, we have that
all morphisms  in this category are epimorphisms. One shows that 
the opposite category $\mathscr{C}_{2}^{\mathrm{op}}$ satisfies
Condition~\ref{pbm}.   For $X\in \mathscr{C}_{1}$ we have that 
$\mathrm{End}_{\mathscr{C}_{1}}\big(X\big)$ is the symmetric group on $X$. 
A direct computation shows that 
$\mathrm{End}_{\mathscr{P}(\mathscr{C}_{2}^{\mathrm{op}})}\big(X\big)$ is the 
dual symmetric inverse monoid $\mathscr{I}^*_X$, and
$\mathrm{End}_{\mathscr{Q}(\mathscr{C}_{2}^{\mathrm{op}})}\big(X\big)$
is the greatest factorizable submonoid $\mathscr{F}^*_X$ of 
$\mathscr{I}^*_X$, see \cite{FL}.

\subsection{The semigroup $\mathcal{PT}_n$ of all partial 
transformations}\label{s3.2}

Let $\mathscr{C}_{3}$ denote the category, whose objects are 
all finite sets and morphisms are all maps between these sets.
One shows that $\mathscr{C}_{3}$ satisfies Condition~\ref{pbm}.
For $X\in \mathscr{C}_{3}$ we have that
$\mathrm{End}_{\mathscr{C}_{3}}\big(X\big)=\mathcal{T}(X)$ 
is the full transformation semigroup on $X$. A direct calculation 
shows that $\mathrm{End}_{\mathscr{P}(\mathscr{C}_{3})}\big(X\big)$
is the semigroup $\mathcal{PT}(X)$ of all partial transformations
on $X$. One can consider the dual situation (as in 
Subsection~\ref{s3.15}). It might be interesting to try to 
understand the resulting monoid.

%

\subsection{Partial endomorphisms of a group}\label{s3.4}

Consider the category $\mathscr{C}_{4}$, whose objects are finite
groups and morphisms are group homomorphisms.  One shows that
$\mathscr{C}_{4}$ satisfies Condition~\ref{pbm}.
If $G$ is a finite group, the monoid 
$\mathrm{End}_{\mathscr{C}_{4}}\big(G\big)$
is the monoid of all endomorphisms of $G$. The monoid
$\mathrm{End}_{\mathscr{P}(\mathscr{C}_{4})}\big(G\big)$
the monoid of all partial endomorphisms of $G$, see \cite{NN}.
Note that $\mathrm{End}_{\mathscr{P}(\mathscr{C}_{4})}\big(G\big)\neq
\mathrm{End}_{\mathscr{Q}(\mathscr{C}_{4})}\big(G\big)$ in general.

\subsection{Partial automorphisms of a group}\label{3.5}

Consider the category $\mathscr{C}_{5}$, whose objects are finite
groups and morphisms are group monomorphisms. One shows that
$\mathscr{C}_{5}$ satisfies Condition~\ref{pbm}.
If $G$ is a finite group,
$\mathrm{End}_{\mathscr{C}_{5}}\big(G\big)$
is the group of all automorphisms of $G$. The monoid 
$\mathrm{End}_{\mathscr{P}(\mathscr{C}_{5})}\big(G\big)$ is the
inverse monoid of all partial automorphisms of $G$. 
The dual construction (using epimorphisms) gives the monoid of 
bicongruences of a group, defined in \cite{Fi} (even in a more 
general setup of  universal algebras).

\subsection{Partial linear endomorphisms}\label{s3.6}

Let $\Bbbk$ be a field. Denote by $\mathscr{C}_{6}$ the category,
whose objects are $\Bbbk^n$, $n=0,1,\dots$, and morphisms
are all linear maps. One shows that $\mathscr{C}_{6}$ satisfies
Condition~\ref{pbm}. The monoid 
$\mathrm{End}_{\mathscr{C}_{6}}\big(\Bbbk^n\big)$ is the
monoid of all $n\times n$ matrices over $\Bbbk$. The monoid
$\mathrm{End}_{\mathscr{P}(\mathscr{C}_{6})}\big(\Bbbk^n\big)$
is the monoid of all partial linear maps on $\Bbbk^n$. 

\subsection{Partial linear automorphisms}\label{s3.7}

Denote by $\mathscr{C}_{7}$ the category,
whose objects are $\Bbbk^n$, $n=0,1,\dots$, and morphisms
are all injective linear maps. One shows that $\mathscr{C}_{7}$ 
satisfies Condition~\ref{pbm}. The monoid 
$\mathrm{End}_{\mathscr{C}_{7}}\big(\Bbbk^n\big)$ is the group
of all invertible $n\times n$ matrices over $\Bbbk$. 
The inverse monoid 
$\mathrm{End}_{\mathscr{P}(\mathscr{C}_{7})}\big(\Bbbk^n\big)$
is the monoid of all partial linear automorphisms of $\Bbbk^n$,
studied in e.g. \cite{Ku}. 

\subsection{Partially defined co-finite automorphisms of 
integers}\label{s3.8}

Consider $\mathbb{Z}$ as a metric space with respect to the
usual metric $\rho(m,n)=|m-n|$. Let $\mathscr{C}_{8}$ be the 
category, whose objects are all co-finite subsets of $\mathbb{Z}$
and morphisms are all isometries. One shows that $\mathscr{C}_{8}$ 
satisfies Condition~\ref{pbm}. The monoid 
$\mathrm{End}_{\mathscr{C}_{8}}\big(\mathbb{Z}\big)$ is the group
of all isometries of $\mathbb{Z}$, which is in fact isomorphic to the
infinite dihedral group. The inverse monoid 
$\mathrm{End}_{\mathscr{P}(\mathscr{C}_{8})}\big(\mathbb{Z}\big)$
is the inverse monoid of all partially defined co-finite
automorphisms of integers, studied in \cite{Be}.

\subsection{Leech's approach to inverse monoids}\label{s3.9}

Let $\mathscr{C}_{9}$ be an {\em abstract division category}
with {\em initial object} $I$ in the sense of \cite[1.1]{Le}. 
Then $\mathscr{C}_{9}^{\mathrm{op}}$ satisfies Condition~\ref{pbm}.
The monoid $\mathrm{End}_{\mathscr{P}(\mathscr{C}_{9}^{\mathrm{op}})}
\big(I\big)$ is isomorphic to the inverse monoid associated with
$\mathscr{C}_{9}$ as defined in \cite[2.1]{Le}. In particular,
if $\mathscr{C}_{9}$ has one object and all morphisms in 
$\mathscr{C}_{9}$ are monomorphisms, then Condition~\ref{pbm}
describes exactly the situation, dual to the one
discussed in \cite[1.2]{Le}.

\subsection{Bisimple inverse monoid}\label{s3.15new}

Let $\mathscr{C}_{10}$ denote the category, whose objects are
cofinite subsets of $\mathbb{N}$ and morphisms are all possible
injections with cofinite image (in particular, all objects of
$\mathscr{C}_{10}$ are isomorphic). In the same way as for 
$\mathscr{C}_{1}$ one shows that $\mathscr{C}_{10}$ satisfies
Condition~\ref{pbm}.  The monoid
$\mathrm{End}_{\mathscr{P}(\mathscr{C}_{10}^{\mathrm{op}})}
\big(\mathbb{N}\big)$ is a bisimple inverse monoid.

\section{The inverse braid monoid}\label{s4}

In this section we give one more example. The semigroup we will
talk about is the {\em inverse braid monoid}, recently defined in
\cite{EL}. But our approach to this monoid will be quite different. 
We will obtain this monoid as one more application of the functor 
$\mathscr{P}$. 

Consider the category $\mathscr{B}$ defined as follows: The objects
of $\mathscr{B}$ are indexed by finite subsets of $\mathbb{N}$.
If $X\subset \mathbb{N}$, then the object, associated with 
$X$, is the plane $P_X=\mathbb{R}^2$ with marked points 
$(i,0)$, $i\in X$.
For $X,Y\subset \mathbb{N}$, $|X|,|Y|\leq \infty$, the 
set $\mathscr{B}(P_X,P_Y)$ is the set of all isotopy classes of 
homeomorphisms with compact support from $P_X$ to $P_Y$, which map
marked points to marked points. The composition is induced by
the usual composition of maps. In particular, if $|X|>|Y|$ then
$\mathscr{B}(P_X,P_Y)=\varnothing$; if $|X|=|Y|=n$ then
$\mathscr{B}(P_X,P_Y)$ is identified (elementwise) with Artin's braid 
group  $\mathcal{B}_n$, see \cite[Theorem~1.10]{Bi}; if $|X|<|Y|$ we have 
$\binom{|Y|}{|X|}$ ways to choose the set $A$ of values for the marked 
points from $P_X$ and after fixing it the part of 
$\mathscr{B}(P_X,P_Y)$ corresponding to $A$ is again
identified (elementwise) with $\mathcal{B}_{|X|}$ via \cite[Theorem~1.10]{Bi}.
We also have that  $\mathscr{B}(P_X,P_X)=\mathrm{B}_{|X|}$ as a group.

First we note that obviously all morphisms in $\mathscr{B}$ are
monomorphisms. Further, $\mathscr{B}$ satisfies
Condition~\ref{pbm}. Indeed, let $f\in \mathscr{B}(P_X,P_Y)$
and $g\in \mathscr{B}(P_Z,P_Y)$. We have $|X|\leq |Y|$ and $|Z|\leq |Y|$.
Let $A$ and $B$ denote the set of marked points in $P_Y$, which
are images of mark points from $P_X$ and $P_Z$ under $f$ and
$g$ respectively. Let $C=A\cap B$, $l=|C|$, and $V\subset \mathbb{N}$
be such that $|V|=l$. Let 
$h\in \mathscr{B}(P_V,P_X)$ be any map, which sends the marked
points to $f^{-1}(C)$. Consider the set $\mathbb{S}$ of all morphisms
from $\mathscr{B}(P_V,P_Z)$, which send marked points
to $g^{-1}(C)$. Then $g\mathbb{S}$ is identified with $\mathcal{B}_l$ 
via \cite[Theorem~1.10]{Bi} and hence there exists a unique
$h'\in X$ such that $gh'=fh$. One now easily shows that
the diagram 
\begin{displaymath}
\xymatrix{
P_X\ar@{^{(}->}[r]^f & P_Y\\
P_V\ar@{^{(}->}[r]^{h'}\ar@{^{(}->}[u]^h & P_Z\ar@{^{(}->}[u]^g 
}
\end{displaymath}
is a pullback. This implies Condition~\ref{pbm}. In particular,
we have the partialization $\mathscr{P}(\mathscr{B})$ and from
Proposition~\ref{pmonoinv} we get that the monoid
$\mathrm{End}_{\mathscr{P}(\mathscr{B})}(P_X)$ is an
inverse monoid.

Let $X=\{1,2,\dots,n\}$. We now claim that
$\mathrm{End}_{\mathscr{P}(\mathscr{B})}(P_X)$ is isomorphic to the 
inverse braid monoid  $\mathcal{IB}_n$ as defined in \cite{EL}. The 
monoid  $\mathcal{IB}_n$ is defined as the set of geometrical braids 
with  $n$ strands in which some strands can be missing with the obvious
multiplication induced by the  multiplication of geometrical braids. 
Let us construct a bijection from $\mathcal{IB}_n$ to
$\mathrm{End}_{\mathscr{P}(\mathscr{B})}(P_X)$ as follows: Denote
the base points for geometric braids by $1,2,\dots,n$. Take some 
partial geometrical braid  $\alpha\in \mathcal{IB}_n$. Assume that 
$\alpha$  consists of  $m$ strands. Then $\alpha$ is given by 
two subsets $A,B\subset \{1,2,\dots,n\}$ and a usual braid 
$\overline{\alpha}$ on $m$ strands, namely, $\overline{\alpha}$ 
connects initial points from $A$ with the terminal points from
in $B$.

The class of the identity transformation of $\mathbb{R}^2$ is 
a morphism in $\mathscr{B}(P_A,P_X)$, which we denote by 
$i(A,X)$. Let $g\in \mathscr{B}(P_A,P_X)$ be the isotopy
class of maps, which maps $\{(j,0)\,:\,j\in A\}$, to 
$\{(k,0)\,:\,k\in B\}$ and corresponds to $\overline{\alpha}$
under the identification, given by \cite[Theorem~1.10]{Bi}.
Then it follows immediately that the map
\begin{eqnarray*}
\mathcal{IB}_n&\rightarrow &
\mathrm{End}_{\mathscr{P}(\mathscr{B})}(P_X) \\
\alpha&\mapsto&
\xymatrix{
P_X & P_A\ar@{_{(}->}[l]_{i(A,X)}\ar[r]^g & P_X
}
\end{eqnarray*}
is in fact an isomorphism of monoids. 

\section{The inverse braid-permutation monoid}\label{s5}

In this section we generalize the example from the previous
section and construct a new inverse monoid with topological origin. 
Again our construction is an immediate application of the functor 
$\mathscr{P}$.

\subsection{Definition}\label{s5.1}

Consider the category $\overline{\mathscr{B}}$ defined as follows: 
The objects of $\overline{\mathscr{B}}$ are indexed by finite subsets 
of $\mathbb{N}$. If $X\subset \mathbb{N}$, then the object, associated 
with  $X$, is the space $Q_X=\mathbb{R}^3$ with marked circles 
$\{(\cos(a),\sin(a),i)\,:\,a\in\left[0,2\pi\right.)\}$, $i\in X$, with the 
orientation induced by the natural order on $\left[0,2\pi)\right.$.
For $X,Y\subset \mathbb{N}$, $|X|,|Y|\leq \infty$, the  set 
$\overline{\mathscr{B}}(Q_X,Q_Y)$ is the set of all isotopy classes 
of  diffeomorphisms with compact support from $Q_X$ to $Q_Y$, which map
marked circles to marked circles preserving the orientation. 
The composition is induced by the usual composition of maps.
In particular, if $|X|>|Y|$ then 
$\overline{\mathscr{B}}(Q_X,Q_Y)=\varnothing$; if $|X|=|Y|=n$ then
$\overline{\mathscr{B}}(Q_X,Q_Y)$ is identified (elementwise) with 
the braid-permutation group  $\mathrm{BP}_n$, see \cite{Da,Wa}
(see also \cite{FRR,Ru}); if $|X|<|Y|$ we have 
$\binom{|Y|}{|X|}$ ways to choose the set $A$ of marked circles in
$Q_Y$ for the values of the marked circles from $Q_X$ and after 
fixing it the part of $\overline{\mathscr{B}}(Q_X,Q_Y)$ 
corresponding to $A$ is again identified (elementwise) with 
$\mathrm{BP}_n$. It follows that
$\overline{\mathscr{B}}(Q_X,Q_X)=\mathrm{BP}_{|X|}$ as a group.

Again all morphisms in $\overline{\mathscr{B}}$ are monomorphisms.
In the same way as in the previous section one shows that 
$\overline{\mathscr{B}}$ satisfies Condition~\ref{pbm}. Hence
we can apply $\mathscr{P}$ and by Proposition~\ref{pmonoinv} we 
get that the monoid
$\mathrm{End}_{\mathscr{P}(\overline{\mathscr{B}})}(Q_X)$ is an
inverse monoid. We call this monoid the 
{\em inverse braid-permutation monoid} and denote it by
$\mathcal{IBP}_n$, where $n=|X|$.

\subsection{Idempotents and factorizability of 
$\mathcal{IBP}_n$}\label{s5.2}

Let $Y\subset X$. Then the class of the identity transformation
of $\mathbb{R}^3$ is a 
monomorphism in $\overline{\mathscr{B}}(Q_Y,Q_X)$, which 
we denote by $f_Y$. Denote by $\varepsilon_Y$
the element $\xymatrix{Q_X & Q_Y\ar@{_{(}->}[l]_{f_Y}
\ar@{^{(}->}[r]^{f_Y} & Q_X}$ in $\mathcal{IBP}_n$.

\begin{lemma}\label{lem61}
$E(\mathcal{IBP}_n)=\{\varepsilon_Y\,:\,Y\subset X\}$.
In particular, $E(\mathcal{IBP}_n)$ is canonically isomorphic 
to the Boolean $(2^X,\cap)$ of $X$. 
\end{lemma}

\begin{proof}
Consider the element $x\in \mathcal{IBP}_n$ given by
$\xymatrix{Q_X & Q_Y\ar@{_{(}->}[l]_{f}
\ar@{^{(}->}[r]^{f} & Q_X}$, where 
$f\in \overline{\mathscr{B}}(Q_Y,Q_X)$. Let $\mathtt{f}\in f$
and $Z=\mathtt{f}(Y)$. Then the class of $\mathtt{f}^{-1}$
is an element of $\overline{\mathscr{B}}(Q_Z,Q_Y)$, call it $g$.
Note that $g$ is an isomorphism since $|Y|=|Z|$.
It follows that the following diagram commutes:
\begin{displaymath}
\xymatrix{Q_X\ar@{=}[d] & Q_Y\ar@/^/[d]^{g^{-1}}\ar@{_{(}->}[l]_{f}
\ar@{^{(}->}[r]^{f} & Q_X\ar@{=}[d]\\
Q_X & Q_Z\ar@{_{(}->}[l]_{f_Z}\ar@/^/[u]^{g}
\ar@{^{(}->}[r]^{f_Z} & Q_X
}
\end{displaymath}
Hence Lemma~\ref{lem:idempotents} implies that 
$E(\mathcal{IBP}_n)=\{\varepsilon_Y\,:\,Y\subset X\}$.
A direct calculation shows that $\varepsilon_Y\mapsto Y$
is an epimorphism from  $E(\mathcal{IBP}_n)$
to $(2^X,\cap)$. Hence it is an isomorphism
since $|E(\mathcal{IBP}_n)|\leq 2^{n}$ by above.
\end{proof}

By construction, we have 
\begin{displaymath}
G(\mathcal{IBP}_n)=\mathfrak{i}(
\mathrm{End}_{\overline{\mathscr{B}}}(Q_X))\cong
\mathrm{BP}_n.
\end{displaymath}

\begin{lemma}\label{lem62}
$\mathcal{IBP}_n$ is a factorizable monoid.
\end{lemma}

\begin{proof}
Let $x\in \mathcal{IBP}_n$. By the same arguments as in
the proof of Lemma~\ref{lem61} we can assume that $x$
is given by $\xymatrix{Q_X & Q_Y\ar@{_{(}->}[l]_{f_Y}
\ar@{^{(}->}[r]^{f} & Q_X}$ for some 
$f\in \overline{\mathscr{B}}(Q_Y,Q_X)$. Choose a
representative $\mathtt{f}\in f$, which preserves the
set of marked circles, given by $X$. Then we can also
consider the class $g\in \overline{\mathscr{B}}(Q_X,Q_X)$ 
containing $\mathtt{f}$. The class $g$ defines an element of
$\mathrm{BP}_n=G(\mathcal{IBP}_n)$, which we call $y$. 
A direct calculation shows that $x=y\varepsilon_Y$. This 
completes the proof.
\end{proof}

\subsection{Presentation of $\mathcal{IBP}_n$}\label{s5.3}

A presentation for  $\mathrm{BP}_n$ was obtained in \cite{Mc,FRR}.
A presentation for $\mathcal{IB}_n$ was obtained in \cite{EL}.
In this subsection we obtain a presentation for $\mathcal{IBP}_n$.
In some sense it is a unification of the results of 
from \cite{FRR} and \cite{EL}. Our arguments are based on the
approach to presentations of factorizable inverse monoids, worked
out in \cite{EEF}. To formulate the result we will need to recall
the presentation of $\mathrm{BP}_n$, obtained in \cite{FRR}.

To simplify notation we assume $X=\{1,2,\dots,n\}$. Let 
$i\in\{1,2,\dots,n-1\}$. Denote by $\tau_i$ the element of
$\mathrm{BP}_n$ given by a homeomorphism of $\mathbb{R}^3$,
whose support contains only the $i$-th and the $(i+1)$-st circles,  
and which interchanges these two circles without moving them 
through each other. Denote by $\sigma_i$ the element of 
$\mathrm{BP}_n$ given by a  homeomorphism of $\mathbb{R}^3$, 
whose support contains only the $i$-th and the $(i+1)$-st circles,  
and which interchanges these two circles by moving the $i$-th
circle through the $(i+1)$-st (see \cite[Section~2]{Ru}). 
It turns out that $\{\sigma_i,\tau_i\,:\,i=1,\dots,n-1\}$
is a generating set for $\mathrm{BP}_n$. We even have:

\begin{theorem}(\cite{FRR})\label{thmfrr}
The group $\mathrm{BP}_n$ is generated by
$\{\sigma_i,\tau_i\,:\,i=1,\dots,n-1\}$ with the following
defining relations:
\begin{equation}\label{eqrel1}
\text{(braid group relations)}\quad
\left\{
\begin{array}{lcl}
\sigma_i\sigma_j & = & \sigma_j\sigma_i,\,\,\,\,\,|i-j|>1;\\
\sigma_i\sigma_{i+1}\sigma_i & = & \sigma_{i+1}\sigma_i\sigma_{i+1};
\end{array}
\right.
\end{equation}
\begin{equation}\label{eqrel2}
\text{(permutation group relations)}\quad
\left\{
\begin{array}{lcl}
\tau_i^2 & = & 1;\\
\tau_i\tau_j & = & \tau_j\tau_i,\,\,\,\,\,|i-j|>1;\\
\tau_i\tau_{i+1}\tau_i & = & \tau_{i+1}\tau_i\tau_{i+1};
\end{array}
\right.
\end{equation}
\begin{equation}\label{eqrel3}
\text{(mixed relations)}\quad
\left\{
\begin{array}{lcl}
\sigma_i\tau_j & = & \tau_j\sigma_i,\,\,\,\,\,|i-j|>1;\\
\tau_i\tau_{i+1}\sigma_i & = & \sigma_{i+1}\tau_i\tau_{i+1};\\
\sigma_i\sigma_{i+1}\tau_i & = & \tau_{i+1}\sigma_i\sigma_{i+1}.
\end{array}
\right.
\end{equation}
\end{theorem}

For $i=1,\dots,n$ denote by $\varepsilon_i$ the element
$\varepsilon_Y$, where $Y=X\setminus\{i\}$. We have the
following:

\begin{theorem}\label{thm65}
The monoid $\mathcal{IBP}_n$ is generated by the elements
$\{\sigma_i,\sigma_i^{-1},\tau_i\,:\,i=1,\dots,n-1\}$ and
$\{\varepsilon_i\,:\,i=1,\dots,n\}$ with the defining 
relations \eqref{eqrel1}--\eqref{eqrel3} and the
following additional relations:
\begin{equation}\label{eqrel305}
\text{(inverse relation)}\quad
\left\{
\begin{array}{lclcl}
\sigma_i\sigma_i^{-1}&=&\sigma_i^{-1}\sigma_i&=&1;
\end{array}
\right.
\end{equation}
\begin{equation}\label{eqrel4}
\text{(semilattice relations)}\quad
\left\{
\begin{array}{lcl}
\varepsilon_i^2&=&\varepsilon_i;\\
\varepsilon_i\varepsilon_j&=&\varepsilon_j\varepsilon_i;
\end{array}
\right.
\end{equation}
\begin{equation}\label{eqrel5}
\text{(action relations)}\quad
\left\{
\begin{array}{lcl}
\sigma_i\varepsilon_j&=&\varepsilon_j\sigma_i,\,\,\,\,\,i\neq j,j+1;\\
\tau_i\varepsilon_j&=&\varepsilon_j\tau_i,\,\,\,\,\,i\neq j,j+1;\\
\sigma_i\varepsilon_i&=&\varepsilon_{i+1}\sigma_i;\\
\sigma_i\varepsilon_{i+1}&=&\varepsilon_{i}\sigma_i;\\
\tau_i\varepsilon_{i}&=&\varepsilon_{i+1}\tau_i;\\
\end{array}
\right.
\end{equation}
\begin{equation}\label{eqrel6}
\text{(kernel relations)}\quad
\left\{
\begin{array}{lclcl}
\varepsilon_{i}\varepsilon_{i+1}\sigma_i&=&
\varepsilon_{i}\varepsilon_{i+1}\tau_i&=&
\varepsilon_{i}\varepsilon_{i+1};\\
\varepsilon_{i}\sigma_i^2&=&\varepsilon_{i}\sigma_i\tau_i&=&
\varepsilon_{i}
\end{array}
\right.
\end{equation}
\end{theorem}

\begin{proof}
The generators
$\{\sigma_i,\sigma_i^{-1},\tau_i\,:\,i=1,\dots,n-1\}$ and
relations \eqref{eqrel1}--\eqref{eqrel305} give a presentation 
for $\mathrm{BP}_n$ by Theorem~\ref{thmfrr}. The relations
\eqref{eqrel4} give a presentation for the semilattice
$(\mathcal{B}(X),\cap)$. The relations \eqref{eqrel5},
which are verified by a direct calculation, give an action of
$\mathrm{BP}_n$ on $(\mathcal{B}(X),\cap)$. 

According to \cite[Theorem~6]{EEF}, the only relations, which are
left are the ones of the form $eg=e$, where 
$e\in E(\mathcal{IBP}_n)$ and $g\in \mathrm{BP}_n$. So, to complete
the proof we have only to show that all these relations can be
derived from the ones described above together with \eqref{eqrel6}
(the latter ones are again easily verified by a direct calculation).

For  $Y\subset X$ set $G_{Y}=\{g\in \mathrm{BP}_n\,:\,
e_Y g=e_Y\}$. By definition, the group $G_{Y}$ consists of classes
of diffeomorphisms with compact support from $Q_X$ to $Q_X$, whose
support does not intersect any of the circles, corresponding to $Y$.
One easily shows that $G_{Y}$ is generated by:
\begin{enumerate}[(a)]
\item\label{ena1} The canonical image of 
$\overline{\mathscr{B}}(Q_{X\setminus Y},Q_{X\setminus Y})$ 
in $\mathrm{BP}_n$, given by
the natural inclusion.
\item\label{ena2} Extra moves, which can  be described as 
follows: Take two circles, the first one for $i\not\in Y$ and the second 
one for $j\in Y$. Move the first circle through the second one and return
it back without involving any other circles. 
\end{enumerate}

The group $\overline{\mathscr{B}}(Q_{X\setminus Y},Q_{X\setminus Y})$ 
is isomorphic to $\mathrm{BP}_m$, where $m=|X|-|Y|$ and hence is 
generated by the  corresponding $\sigma_i$'s and $\tau_i$'s. Hence 
in the equality $eg=e$ we can assume $g$ to be one of these generators. 
Using conjugation by elements with $\mathrm{BP}_n$ (i.e. relations
\eqref{eqrel5}) we thus reduce $eg=e$ to the first line of the
kernel relations \eqref{eqrel6}.

The extra moves described in \eqref{ena2} can be of two different 
kinds: after moving the first circle through the second one we can 
return it back with or without moving the second circle 
through the first one. This means that these extra moves are
either of the form $x\sigma_i^{2} x^{-1}$ or of the form
$x\sigma_i\tau_i x^{-1}$ (or  $x\tau_i\sigma_i x^{-1}$), 
where $x\in \mathrm{BP}_n$. Note that 
$\varepsilon_{i}\tau_i\sigma_i=\varepsilon_{i}$ follows
from $\varepsilon_{i}\sigma_i\tau_i=\varepsilon_{i}$,
$\tau_i^2=1$ and $\varepsilon_{i}\sigma_i^2=\varepsilon_{i}$.
Hence, up to
conjugation by elements with $\mathrm{BP}_n$ (i.e. relations
\eqref{eqrel5}) the condition $eg=e$, where $g$ is our extra move,
reduces to the second line of the kernel relations \eqref{eqrel6}.

Now the statement of Theorem~\ref{thm65} follows immediately from
\cite[Theorem~6]{EEF}.
\end{proof}

\begin{remark}\label{rem68}
{\rm
The system of generators of $\mathcal{IBP}_n$, presented in
Theorem~\ref{thm65} is reducible: relations \eqref{eqrel5} show
that $\mathcal{IBP}_n$ is already generated by
$\{\sigma_i,\sigma_i^{-1},\tau_i\,:\,i=1,\dots,n-1\}$ and
the element $\varepsilon=\varepsilon_1$. 
From Theorem~\ref{thm65} one easily derives
that with repsect to this irreducible system of generators
the defining relations are \eqref{eqrel1}-\eqref{eqrel305} 
together with the following additional relations:
\begin{displaymath}
\left\{
\begin{array}{lclclclcll}
\varepsilon^2&=&\varepsilon&=&\varepsilon\sigma_1^2&=&
\sigma_1^2\varepsilon&=&\varepsilon\sigma_1\tau_1;\\
\varepsilon\sigma_i&=&\sigma_i\varepsilon,&&i>2;\\
\varepsilon\tau_i&=&\tau_i\varepsilon,&&i>2;\\
\varepsilon\sigma_1\varepsilon&=&
\varepsilon\sigma_1\varepsilon\sigma_1&=&
\sigma_1\varepsilon\sigma_1\varepsilon.
\end{array}
\right.
\end{displaymath}
}
\end{remark}

\begin{remark}\label{rem69}
{\rm
The same approach as we used in Theorem~\ref{thm65} can be
used to derive a presentation for $\mathcal{IB}_n$,
substantially shortening the arguments from \cite{EL}.
}
\end{remark}

\section{Iterations of  the partialization functor}\label{s6}

\subsection{$\mathscr{P}$ and monomorphisms}\label{s6.1}

To be able to iterate $\mathscr{P}$ (or $\mathscr{Q}$) one has
to ensure that $\mathscr{P}(\mathscr{C})$ (or 
$\mathscr{Q}(\mathscr{C})$ respectively) satisfies
Condition~\ref{pbm}. This is wrong in the general case, we
will give an example in Subsection~\ref{s6.15}. The first 
important step to understand Condition~\ref{pbm} for 
$\mathscr{P}(\mathscr{C})$  ($\mathscr{Q}(\mathscr{C})$) is 
the following result, which describes monomorphisms in
partialized categories:

\begin{proposition}\label{p12}
\begin{enumerate}[(a)]
\item\label{p12.1}
Let $x=\mathrm{D}(\mi,\mj,\mk,\alpha,f)$. 
Then $x$ is a monomorphism if and only if 
$x=\mathfrak{i}(\gamma)$ for some monomorphism
$\gamma\in \mathscr{C}(\mi,\mj)$.
\item\label{p12.2}
Let $y=\overline{\mathrm{D}}(\mi,\mj,\mk,\alpha,f,g)$. 
Then $y$ is a monomorphism if and only if 
$x=\mathfrak{i}'(\gamma)$ for some monomorphism
$\gamma\in \mathscr{C}(\mi,\mj)$.
\end{enumerate}
\end{proposition}

\begin{proof}
The statement \eqref{p12.2} is a special case of the statement
\eqref{p12.1}, so we prove \eqref{p12.1}. Suppose $x$ is a 
monomorphism. Let $a=\mathrm{D}(\mi,\mi,\mk,\alpha,\alpha)$ and
$b=\mathrm{D}(\mi,\mi,\mi,\mathtt{1}_{\mi},\mathtt{1}_{\mi})$. 
By a direct  calculation one obtains  $xa=xb=x$, implying that $a=b$. In
particular, $\alpha$ is an isomorphism. Without loss of generality 
we hence can assume $x=\mathfrak{i}(f)$.
Let $g_1,g_2\in \mathscr{C}(\ml,\mi)$ be such that 
$g_1\neq g_2$. Since $\mathfrak{i}$ is injective, we have
$\mathfrak{i}(g_1)\neq \mathfrak{i}(g_1)$. Since
$x$ is a monomorphism  in $\mathscr{P}(\mathscr{C})$ we
get  $\mathfrak{i}(fg_1)=x\mathfrak{i}(g_1)
\neq x\mathfrak{i}(g_1)=\mathfrak{i}(fg_2)$.
Hence $fg_1\neq fg_2$. This implies that $f$ is 
a monomorphism in $\mathscr{C}$.

Now let $f\in \mathscr{C}(\mi,\mj)$ be a monomorphism
and $x=\mathfrak{i}(f)$.
Let $c=\mathrm{D}(\mk,\mi,\ml,\alpha,g)$ and
$d=\mathrm{D}(\mk,\mi,\ml',\beta,h)$ be such that
$xc=xd$. Hence there exists an isomorphism
$\gamma\in \mathscr{C}(\ml,\ml')$ such that
the solid part of the following diagram commutes:
\begin{displaymath}
\xymatrix{
\mk\ar@{=}[d]&\ml\ar[r]^{g}\ar@{_{(}->}[l]_{\alpha}
\ar@/^/[d]^{\gamma}
&\mi\ar@{^{(}->}[r]^{f}\ar@{::}[d]&\mj\ar@{=}[d]\\
\mk&\ml'\ar[r]_{h}\ar@{_{(}->}[l]^{\beta}
\ar@/^/[u]^{\gamma^{-1}}
&\mi\ar@{^{(}->}[r]_{f}&\mi\\
}
\end{displaymath}
Since $f$ is a monomorphism, the middle square of the 
diagram commutes as well, implying that $c$ and $d$
are equivalent. This means that $x$ is a monomorphism.
\end{proof}

\subsection{Quasi-iterations of $\mathscr{P}$}\label{s6.01}

Let $\mathscr{C}$ be a category satisfying Condition~\ref{pbm}.
Then the category $\mathscr{P}(\mathscr{C})$ is well-defined,
however, it does not have to satisfy Condition~\ref{pbm} (this
problem will be addressed later on in this section). Anyway,
we can just formally consider the morphisms in the 
``second partialization'' $\mathscr{P}(\mathscr{P}(\mathscr{C}))$
as defined in \eqref{eqeq1}. These will be elements
$\mathrm{D}(\mathtt{i},\mathtt{j},\mathtt{k},\alpha,f)$, where
$f$ is a morphism from $\mathscr{P}(\mathscr{C})$, hence, in turn,
have the form $\mathrm{D}(\mathtt{k},\mathtt{j},\mathtt{l},\beta,g)$
for some morphism $g$ from $\mathscr{C}$. Because of Proposition~\ref{p12}, 
the element $\mathrm{D}(\mathtt{i},\mathtt{j},\mathtt{k},\alpha,f)$ can be
viewed as the diagram
\begin{displaymath}
\xymatrix{
\mathtt{i} & \mathtt{k}\ar@{_{(}->}[l]_{\alpha}& 
\mathtt{l}\ar[r]^g\ar@{_{(}->}[l]_{\beta} & \mathtt{j},
}
\end{displaymath}
which consists only of elements from the original category 
$\mathscr{C}$. This suggests the following definition:
Let $n\in\mathbb{N}$. We define the category
$\mathscr{P}^{(k)}(\mathscr{C})$ as follows: The objects of
$\mathscr{P}^{(k)}(\mathscr{C})$ are the same as the
same objects as $\mathscr{C}$. For $\mathtt{i},\mathtt{j}\in
\mathscr{C}$ the set $\mathscr{P}^{(k)}(\mathscr{C})(\mathtt{i},\mathtt{j})$
is the set of equivalence classes of diagrams
\begin{equation}\label{eqel1}
\mathrm{D}(n,\mathtt{i},\mathtt{j},\mathtt{k}_i,\alpha_i,f):\quad\quad
\xymatrix{
\mathtt{i} & \mathtt{k}_n\ar@{_{(}->}[l]_{\alpha_n}&
\dots\ar@{_{(}->}[l]_{\alpha_{n-1}}
& \mathtt{k}_{1}\ar@{_{(}->}[l]_{\alpha_1}\ar[r]^f & \mathtt{j},
}
\end{equation}
where $\mathrm{D}(n,\mathtt{i},\mathtt{j},\mathtt{k}_i,\alpha_i,f)$
and $\mathrm{D}(n,\mathtt{i}',\mathtt{j}',\mathtt{k}'_i,\alpha'_i,f')$
are said to be equivalent if there exist isomorphisms
$\gamma_i$ making the following diagram commutative:
\begin{equation}\label{eqel2}
\xymatrix{
\mathtt{i}\ar@{=}[d] & \mathtt{k}_n\ar@{_{(}->}[l]_{\alpha_n}
\ar@/^/@{.>}[d]|-{\gamma_n}&
\dots\ar@{_{(}->}[l]_{\alpha_{n-1}}
& \mathtt{k}_{1}\ar@{_{(}->}[l]_{\alpha_1}\ar[r]^f
\ar@/^/@{.>}[d]|-{\gamma_1} & 
\mathtt{j}\ar@{=}[d]\\
\mathtt{i}' & \mathtt{k}'_n\ar@{_{(}->}[l]_{\alpha'_n}
\ar@/^/@{.>}[u]|-{\gamma_n^{-1}}&
\dots\ar@{_{(}->}[l]_{\alpha'_{n-1}}
& \mathtt{k}'_{1}\ar@{_{(}->}[l]_{\alpha'_1}
\ar[r]^{f'}\ar@/^/@{.>}[u]|-{\gamma_1^{-1}} & \mathtt{j}'\\
}
\end{equation}
Define the composition of 
$\mathrm{D}(n,\mathtt{i},\mathtt{j},\mathtt{l}_i,\alpha_i,f)$
and 
$\mathrm{D}(n,\mathtt{j},\mathtt{k},\mathtt{m}_i,\beta_i,g)$
as the element
$\mathrm{D}(n,\mathtt{i},\mathtt{k},\mathtt{n}_i,\gamma_i,gf_1)$
where $\gamma_n=\alpha_n\dots\alpha_1\gamma'_n$
from the following diagram, where all squares all pullbacks:
\begin{equation}\label{newcomp}
\xymatrix{
&&&&&&\mathtt{k}\\
&&&\mathtt{j}&
\mathtt{m}_n\ar@{_{(}->}[l]_{\beta_n}&
\dots\ar@{_{(}->}[l]_{\beta_{n-1}}
& \mathtt{m}_{1}\ar@{_{(}->}[l]_{\beta_1}\ar[u]^g 
\\
\mathtt{i} & \mathtt{l}_n\ar@{_{(}->}[l]_{\alpha_n}&
\dots\ar@{_{(}->}[l]_{\alpha_{n-1}}
& \mathtt{l}_{1}\ar@{_{(}->}[l]_{\alpha_1}\ar[u]^f & 
\mathtt{n}_{n}\ar@{.>}[l]_{\gamma'_n}\ar@{.>}[u]^{f_n}
&\dots\ar@{.>}[l]_{\gamma_{n-1}}&
\mathtt{n}_{1}\ar@{.>}[l]_{\gamma_1}\ar@{.>}[u]^{f_1}
}
\end{equation}
It is straightforward to verify that the above composition is
well-defined and associative. In particular,
$\mathscr{P}^{(k)}(\mathscr{C})$ is a category. We have the 
canonical inclusion
$\mathfrak{i}^{(k)}:\mathscr{C}\to \mathscr{P}^{(k)}(\mathscr{C})$,
which is defined by sending $f\in \mathscr{C}(\mathtt{i},\mathtt{j})$
to the element 
$\mathrm{D}(n,\mathtt{i},\mathtt{j},\mathtt{i},\mone_{\mathtt{i}},f)$.

If $\mathscr{C}$ and $\mathscr{C}'$ are two categories, satisfying
Condition~\ref{pbm} and $F:\mathscr{C}\to \mathscr{C}'$ is  a functor,
which preserves monomorphisms, then $F$ defines a functor from
$\mathscr{P}^{(k)}(\mathscr{C})$ to 
$\mathscr{P}^{(k)}(\mathscr{C}')$ by mapping
$\mathrm{D}(n,\mathtt{i},\mathtt{j},\mathtt{l}_i,\alpha_i,f)$
to 
$\mathrm{D}(n,F(\mathtt{i}),F(\mathtt{j}),F(\mathtt{l}_i),
F(\alpha_i),F(f))$. In particular, $\mathscr{P}^{(k)}$ becomes
a functor from the category of all small categories, satisfying 
Condition~\ref{pbm}, where morphisms are all functors, preserving 
monomorphisms, to the category of all small categories. We will
call $\mathscr{P}^{(k)}$ the {\em $k$-th quasi-iteration} of
$\mathscr{P}^{(k)}$. The functor $\mathfrak{i}^{(k)}$ defined above
is a natural transformation from the canonical inclusion
functor to $\mathscr{P}^{(k)}$. It is convenient to let
$\mathscr{P}^{(0)}$ denote the identity functor (or the natural
inclusion into the category of all small categories).

For $k\geq 0$ we define a natural transformation 
$\mathfrak{j}_k: \mathscr{P}^{(k)}\to \mathscr{P}^{(k+1)}$
via
\begin{displaymath}
\xymatrix{
&&&\mathtt{j}\\
\mathtt{i} & \mathtt{k}_n\ar@{_{(}->}[l]_{\alpha_n}&
\dots\ar@{_{(}->}[l]_{\alpha_{n-1}}
& \mathtt{k}_{1}\ar@{_{(}->}[l]_{\alpha_1}\ar[u]^f
}\mapsto
\xymatrix{
&&&&\mathtt{j}\\
\mathtt{i} & \mathtt{k}_n\ar@{_{(}->}[l]_{\alpha_n}&
\dots\ar@{_{(}->}[l]_{\alpha_{n-1}}
& \mathtt{k}_{1}\ar@{_{(}->}[l]_{\alpha_1} & 
\mathtt{k}_{1}\ar@{=}[l]\ar[u]^f
}
\end{displaymath}
It is straightforward to verify that $\mathfrak{j}_k$
is indeed a natural transformation, moreover, it is
injective. We obviously have $\mathfrak{j}_0=\mathfrak{i}$.

For $k>1$ we define a natural transformation 
$\mathfrak{p}_k: \mathscr{P}^{(k)}\to \mathscr{P}^{(k-1)}$
via
\begin{displaymath}
\xymatrix{
&&&\mathtt{j}\\
\mathtt{i} & \mathtt{k}_n\ar@{_{(}->}[l]_{\alpha_n}&
\dots\ar@{_{(}->}[l]_{\alpha_{n-1}}
& \mathtt{k}_{1}\ar@{_{(}->}[l]_{\alpha_1}\ar[u]^f
}\mapsto
\xymatrix{
&&&&\mathtt{j}\\
\mathtt{i} & \mathtt{k}_n\ar@{_{(}->}[l]_{\alpha_n}&
\dots\ar@{_{(}->}[l]_{\alpha_{n-1}}
& \mathtt{k}_{3}\ar@{_{(}->}[l]_{\alpha_3}
& \mathtt{k}_{1}\ar@{_{(}->}[l]_{\alpha_2\alpha_1}\ar[u]^f
}
\end{displaymath}
It is straightforward to verify that $\mathfrak{p}_k$
is indeed a natural transformation, moreover, it is surjective.

We have the following picture of functors and natural
transformations:
\begin{equation}\label{eq71-1}
\xymatrix{
\mathscr{P}^{(0)}\ar@/_/@{^{(}->}[r]_{\mathfrak{j}_0} &
\mathscr{P}^{(1)}\ar@/_/@{^{(}->}[r]_{\mathfrak{j}_1} &
\mathscr{P}^{(2)}\ar@/_/@{^{(}->}[r]_{\mathfrak{j}_2}
\ar@/_/@{->>}[l]_{\mathfrak{p}_2}&
\mathscr{P}^{(3)}\ar@/_/@{^{(}->}[r]_{\mathfrak{j}_3}
\ar@/_/@{->>}[l]_{\mathfrak{p}_3}&
\dots\ar@/_/@{->>}[l]_{\mathfrak{p}_4}
}
\end{equation}
Furthermore, we have the following easy fact, which follows
directly from the definitions:

\begin{lemma}\label{l71}
For each $k\geq 1$ the composition
$\mathfrak{p}_{k+1}\mathfrak{j}_k$ is the identity.
\end{lemma}

\begin{remark}\label{r72}
{\rm
Although \eqref{eq71-1} and Lemma~\ref{l71} is enough for our purposes,
which will be explained in the next subsections, the real structure of
$\{\mathscr{P}^{(k)}\}$ is richer. In fact, for each $k\geq 1$ and each
$s$, $1\leq s\leq k$, we can define the natural transformation
$\mathfrak{j}_k^{(s)}:\mathscr{P}^{(k)}\to\mathscr{P}^{(k+1)}$
as follows: In the element 
$\mathrm{D}(n,\mathtt{i},\mathtt{j},\mathtt{k}_i,\alpha_i,f)$
we substitute the fragment $\mathtt{k}_s$ with the
fragment $\xymatrix{\mathtt{k}_s\ar@{=}[r]&\mathtt{k}_s}$.
Further, for each $k\geq 2$ and each $s$, $2\leq s\leq k$, we 
can define the natural transformation 
$\mathfrak{p}_k^{(s)}:\mathscr{P}^{(k)}\to\mathscr{P}^{(k-1)}$
as follows: In the element 
$\mathrm{D}(n,\mathtt{i},\mathtt{j},\mathtt{k}_i,\alpha_i,f)$
we skip $\mathtt{k}_s$ and let the map from $\mathtt{k}_{s-1}$ 
to $\mathtt{k}_{s+1}$ be $\alpha_s\alpha_{s-1}$ (here by
$\mathtt{k}_{k+1}$ we mean $\mathtt{i}$). In particular, in the 
above notation we have $\mathfrak{j}_k=\mathfrak{j}_k^{(1)}$ and
$\mathfrak{p}_k=\mathfrak{p}_k^{(2)}$. For all $k$ and all 
appropriate $s$ we have that the compositions
$\mathfrak{p}_{k+1}^{(s)}\mathfrak{j}_k^{(s)}$ and
$\mathfrak{p}_{k+1}^{(s+1)}\mathfrak{j}_k^{(s)}$ are identities.
We will describe in more details the analogous structure for 
$\mathscr{P}^{k}$ later on in Subsection~\ref{s6.3}. In fact 
the analogous structure for $\mathscr{P}^{k}$ is even richer.
}
\end{remark}

\subsection{$\mathscr{P}^{(k)}$ and regular semigroups}\label{s6.03}

\begin{theorem}\label{thm72}
Let $\mathscr{C}$ be a category in which all morphisms are
monomorphisms and which satisfies Condition~\ref{pbm}. 
\begin{enumerate}[(i)]
\item \label{thm72.1}
For each $k\geq 0$ and for each $\mathtt{i}\in \mathscr{C}$ the
monoid $\mathrm{End}_{\mathscr{P}^{(k)}(\mathscr{C})}
\big(\mathtt{i}\big)$ is regular.
\item \label{thm72.2}
For each $k\geq 1$ and for each $\mathtt{i}\in \mathscr{C}$ the
monoid $\mathrm{End}_{\mathscr{P}^{(k)}(\mathscr{C})}
\big(\mathtt{i}\big)$ is a retract of the monoid
$\mathrm{End}_{\mathscr{P}^{(k+1)}(\mathscr{C})}
\big(\mathtt{i}\big)$.
\end{enumerate}
\end{theorem}

\begin{proof}
In the case $k=2$ the statement \eqref{thm72.1} follows from the 
following diagram, in which all squares are pullbacks:
\begin{displaymath}
\xymatrix{
&&&&&&\mathtt{j}\\
&&&&\mathtt{i}&\mathtt{l}\ar@{_{(}->}[l]_{\alpha}&
\mathtt{k}\ar@{^{(}->}[u]^{\gamma}\ar@{_{(}->}[l]_{\beta}\\
&&\mathtt{j}&\mathtt{k}\ar@{_{(}->}[l]_{\gamma}&
\mathtt{k}\ar@{^{(}->}[u]^{\alpha\beta}\ar@{=}[l]&&\\
\mathtt{i}&\mathtt{l}\ar@{_{(}->}[l]_{\alpha}
&\mathtt{k}\ar@{^{(}->}[u]^{\gamma}\ar@{_{(}->}[l]_{\beta}
&\mathtt{k}\ar@{==}[l]\ar@{==}[u]&
\mathtt{k}\ar@{==}[l]\ar@{==}[u]&
\mathtt{k}\ar@{==}[l]\ar@{^{(}-->}[uu]^{\beta}&
\mathtt{k}\ar@{==}[l]\ar@{==}[uu]
}
\end{displaymath}
In the general case the argument is the same (but requires
more space to draw).

The statement \eqref{thm72.1} follows immediately from 
\eqref{eq71-1} and Lemma~\ref{l71}.
\end{proof}

Applying Theorem~\ref{thm72} to the category
$\mathscr{C}_1$ from Subsection~\ref{s3.1} we obtain
a series of regular monoids for which  $\mathcal{IS}_n$ is
a retract (and such that each monoid in the series is a 
retract of the next one). These monoids might be interesting 
objects to study. Later on in Subsection~\ref{s6.4} we shall
discuss slightly different orthodox generalizations of 
$\mathcal{IS}_n$.

\subsection{$\mathscr{P}(\mathscr{C})$ does not have to
satisfy Condition~\ref{pbm}}\label{s6.15}

Here we give an example of a category $\mathscr{C}$
satisfying Condition~\ref{pbm} such that $\mathscr{P}(\mathscr{C})$ 
does not satisfy this condition. Let the objects of $\mathscr{C}$ 
be the set ${\mathbb N}$ and all its finite subsets. Set
\begin{displaymath}
\mathscr{C}(X,Y)=
\begin{cases}
\text{all injections from }X \text{ to }Y,& X\neq {\mathbb N}
\text{ or } Y\neq {\mathbb N}\\
\mone_{{\mathbb N}}, & X=Y={\mathbb N}.
\end{cases}
\end{displaymath}
One easily checks that $\mathscr{C}$ is a category. Exactly in the
same way as in Subsection~\ref{s3.1} one shows that 
$\mathscr{C}$ satisfies Condition~\ref{pbm}. 

Let $A=\{1\}$, $B=\{1,2\}$ and 
$f=\mathrm{D}({\mathbb N},B,B,\mathrm{incl},\mone_B)$. 
We claim that the solid part of the 
following digram in $\mathscr{P}(\mathscr{C})$ does not
have a pullback (note that $\mathrm{incl}$ is a monomorphism in
$\mathscr{P}(\mathscr{C})$ by Proposition \ref{p12}):
\begin{equation}\label{epullb}
\xymatrix{
\mathbb{N}\ar[r]^f& B\\
D\ar@{.>}[r]^g\ar@{.>}[u]^{\alpha}& A \ar@{_{(}->}[u]_{\mathrm{incl}}
}
\end{equation}
Assume that this is not the case. Fix some  $k\in\{1,3,4,\dots\}$.
Set $D_k=\{1,3,\dots,k\}$ and denote by $\iota_k$ the natural
inclusion $D_k\hookrightarrow \mathbb{N}$. 
Let $g=g_k=\mathrm{D}(D_k,A,A,\mathrm{incl},\mone_A)$,
$\alpha=\alpha_k=\mathrm{D}(D_k,\mathbb{N},D_k,\mone_D,\iota_k)$.
A direct calculation shows that the diagram \eqref{epullb} commutes
for each $k\in\{1,3,4,\dots\}$. If a pullback would exist,
one easily checks that for the pullback $D\neq \varnothing, \mathbb{N}$.
Further, one shows that, without loss of generality, in the pullback 
we have $D=D_k$ for some $k\in\{1,3,4,\dots\}$ and even
$g=g_k$ and $\alpha=\alpha_k$. The pullback condition and the 
commutativity of our digram for $k+1$ would now imply the existence 
of a map, $\gamma\in \mathscr{P}(\mathscr{C})(D_{k+1},D_k)$ such that
$\iota_{k+1}=\iota_k\gamma$. A direct computation shows that
such $\gamma$ does not exist. A contradiction. Hence
$\mathscr{P}(\mathscr{C})$  does not satisfy  Condition~\ref{pbm}.
 
\subsection{$\mathscr{P}^n$ as an endofunctor}\label{s6.2}

We would like to define some category on which 
$\mathscr{P}$ would be an endofunctor.  After the previous 
subsection it is clear that we can not just take the category
of all categories, satisfying Condition~\ref{pbm}. Hence we impose
one more condition, which at first glance looks rather artificial.
The naturality of this condition will become clear later on, when
we show that it works.

Assume that $\mathscr{C}$ is a category, satisfying 
Condition~\ref{pbm} and the following condition:

\begin{condition}\label{pbm2}
For each $\mathtt{i},\mathtt{j},\mathtt{k}\in \mathscr{C}$ and
monomorphisms $\alpha\in \mathscr{C}(\mathtt{i},\mathtt{j})$
and $\beta\in \mathscr{C}(\mathtt{j},\mathtt{k})$ there exists
$\mathtt{l}\in \mathscr{C}$ and monomorphisms 
$\gamma\in \mathscr{C}(\mathtt{i},\mathtt{l})$ and
$\delta\in \mathscr{C}(\mathtt{l},\mathtt{k})$ such that:
\begin{enumerate}[(a)]
\item\label{pbm2.1} The solid square on the diagram 
\eqref{eqpbm2} is a pullback:
\begin{equation}\label{eqpbm2}
\xymatrix{
\mathtt{n}\ar@{^{(}-->}[rr]^{\zeta}\ar@{.>}[rd]^{\varphi} 
& & \mathtt{k}\\
 & \mathtt{l}\ar@{^{(}->}[ur]^{\delta}& \\
\mathtt{m}\ar@{^{(}-->}[uu]^{\eta}\ar@{^{(}-->}[r]^{\xi} & \mathtt{i}\ar@{^{(}->}[r]^{\alpha}\ar@{^{(}->}[u]^{\gamma}& 
\mathtt{j}\ar@{_{(}->}[uu]_{\beta}
}
\end{equation}
\item\label{pbm2.2} For every $\mathtt{m},\mathtt{n}\in\mathscr{C}$
and monomorphisms $\xi\in \mathscr{C}(\mathtt{m},\mathtt{i})$,
$\eta\in \mathscr{C}(\mathtt{m},\mathtt{n})$ and
$\zeta\in \mathscr{C}(\mathtt{n},\mathtt{k})$ such that the outer square
on the diagram \eqref{eqpbm2} commutes and is a pullback, there exists
$\varphi\in  \mathscr{C}(\mathtt{n},\mathtt{l})$ making the whole
diagram \eqref{eqpbm2} commutative.
\end{enumerate}
\end{condition}

It is easy to see that in Condition~\ref{pbm2} the map
$\varphi$, if exists,  is automatically unique and a monomorphism.
Moreover, the left square of the diagram is a pullback
(since the outer one is and the diagram commutes).
It is also straightforward to verify that
$\mathtt{l}$, $\gamma$ and $\delta$ are defined uniquely up
to an isomorphism. We will call the right square of
\eqref{eqpbm2} the {\em complement diagram}.
As an example, later on in Subsection~\ref{s6.4} we will show 
that the category $\mathscr{C}_1$ from Subsection~\ref{s3.1}
satisfies Condition~\ref{pbm2}.

Denote by $\mathscr{S}$ the category, whose objects are small 
categories, satisfying Condition~\ref{pbm} and Condition~\ref{pbm2},
and whose morphisms are all possible functors, which
preserve monomorphisms.

\begin{theorem}\label{thm81}
The functor $\mathscr{P}$ is an endofunctor of the
category $\mathscr{S}$.
\end{theorem}

\begin{proof}
Let $\mathscr{C}\in \mathscr{S}$. Then Proposition~\ref{p12}
guarantees that $\mathscr{P}(\mathscr{C})$ satisfies 
Condition~\ref{pbm2}. Let $\mathscr{C}'\in \mathscr{S}$ and
$F:\mathscr{C}\to \mathscr{C}'$ be a functor, which
preserves monomorphisms. Then again Proposition~\ref{p12}
guarantees that $\mathscr{P}(F)$ preserves monomorphisms.
Hence we need only to check that $\mathscr{P}(\mathscr{C})$ 
satisfies  Condition~\ref{pbm}.

Because of Proposition~\ref{p12} we can identify the monomorphisms
in the categories $\mathscr{C}$ and $\mathscr{P}(\mathscr{C})$.
Let $x=\mathrm{D}(\mathtt{i},\mathtt{j},\mathtt{k},\alpha,f)$
and $\beta\in \mathscr{C}(\mathtt{l},\mathtt{j})$ be a monomorphism.
Consider the pullback diagram given by the right square of the
diagram \eqref{eq81-1}:
\begin{equation}\label{eq81-1}
\xymatrix{
\mathtt{i} & \mathtt{k}\ar[r]^f\ar@{_{(}->}[l]_{\alpha} 
& \mathtt{j} \\ 
\mathtt{n}\ar@{^{(}-->}[u]^{\beta''} & \mathtt{m}\ar[r]^{f'}
\ar@{_{(}->}[u]_{\beta'}\ar@{_{(}-->}[l]_{\alpha'} & 
\mathtt{l}\ar@{_{(}->}[u]_{\beta} 
}
\end{equation}
Applying Condition~\ref{pbm2} to the monomorphism $\alpha$
and $\beta'$ we obtain $\mathtt{n}$, $\alpha'$ and $\beta''$.
The upper row of \eqref{eq81-1} is the element $x$.
The lower row of \eqref{eq81-1} is an element from
$\mathscr{P}(\mathscr{C})(\mathtt{n},\mathtt{l})$, say $y$.
Since both squares of \eqref{eq81-1}  are pullbacks, a direct
calculation shows that we have the following commutative diagram
in $\mathscr{P}(\mathscr{C})$:
\begin{equation}\label{eq81-2}
\xymatrix{
\mathtt{i}\ar[r]^x & \mathtt{j} \\
\mathtt{n}\ar@{^{(}-->}[u]^{\beta''}\ar[r]^y & 
\mathtt{l}\ar@{_{(}->}[u]_{\beta}
}
\end{equation}
We claim that \eqref{eq81-2} is in fact a pullback.
Assume that for some $\mathtt{p}\in\mathscr{C}$
there is  a monomorphism $\delta\in\mathscr{C}(\mathtt{p},\mathtt{i})$
and $z\in \mathscr{P}(\mathscr{C})(\mathtt{p},\mathtt{l})$
such that $\beta z=x\delta$. Suppose that
$z=\mathrm{D}(\mathtt{p},\mathtt{l},\mathtt{q},\xi,g)$
and consider the following diagram,
the solid part of which commutes because of \eqref{eq81-1}:
\begin{equation}\label{eq81-3}
\xymatrix{
&&\mathtt{i} && \mathtt{k}\ar[rr]^f\ar@{_{(}->}[ll]_{\alpha} 
&& \mathtt{j} \\ 
\mathtt{p}\ar@{^{(}->}[rru]|-{\delta}\ar@{.>}[rr]|-{\varphi}
&&\mathtt{n}\ar@{^{(}->}[u]|-{\beta''} && \mathtt{m}\ar[rr]^{f'}
\ar@{_{(}->}[u]_{\beta'}\ar@{_{(}->}[ll]|->>>>>{\alpha'} && 
\mathtt{l}\ar@{_{(}->}[u]_{\beta} \\
&&\mathtt{q}\ar@{_{(}->}[llu]|-{\xi}\ar[rrrru]|-g
\ar@{^{(}-->}[rru]|-{\eta}\ar@{^{(}-->}[rruu]^>>>>>>>{\delta'}&&&&
}
\end{equation}
Using \eqref{neqeq2} and \eqref{eqeq2} one shows that 
the condition $\beta z=x\delta$ implies the existence of 
$\delta'$ as on \eqref{eq81-3} such that the part  of the
digram \eqref{eq81-3} formed by all solid arrows and
$\delta'$ commutes and the square $\alpha\delta'=\delta\xi$ is a 
pullback. Now the fact that the right square is a pullback
implies the existence of $\eta$ as on \eqref{eq81-3} such that
the part of \eqref{eq81-3} formed by all solid arrows,
$\delta'$ and $\eta$ commutes. Finally, Condition~\ref{pbm2}
implies now the existence of $\varphi$ as on \eqref{eq81-3}
such that the whole digram \eqref{eq81-3} commutes.

The commutativity of \eqref{eq81-3} implies
$\delta=\beta''\varphi$ and $y\varphi=z$. Hence
the diagram \eqref{eq81-2} is a pullback. This
completes the proof.
\end{proof}

\subsection{Multiplication for $\mathscr{P}^k$}\label{s6.25}

Let $\mathscr{C}\in\mathscr{S}$. Then, by Theorem~\ref{thm81}
we have that $\mathscr{P}^n(\mathscr{C})\in\mathscr{S}$ for
all $n\geq 0$ (we assume $\mathscr{P}^0=\mathrm{ID}$). 
By induction one gets that the morphisms in 
$\mathscr{P}^n(\mathscr{C})$ are exactly the equivalence classes
of the diagrams \eqref{eqel1} with respect to the equivalence
defined on \eqref{eqel2}. However, the multiplication of these
elements is quite different from the multiplication in
$\mathscr{P}^{(n)}(\mathscr{C})$ described in \eqref{newcomp}.
The product in $\mathscr{P}^n(\mathscr{C})$ in terms of the
original catgory $\mathscr{C}$ is described in the following
statement:

\begin{proposition}\label{prop91}
Let $\mathrm{D}(n,\mathtt{i},\mathtt{j},\mathtt{l}_i,\alpha_i,f)$
and $\mathrm{D}(n,\mathtt{j},\mathtt{k},\mathtt{m}_i,\beta_i,g)$
be two morphisms in $\mathscr{P}^n(\mathscr{C})$. The product
$\mathrm{D}(n,\mathtt{j},\mathtt{k},\mathtt{m}_i,\beta_i,g)
\mathrm{D}(n,\mathtt{i},\mathtt{j},\mathtt{l}_i,\alpha_i,f)$
in $\mathscr{P}^n(\mathscr{C})$ is the element
$\mathrm{D}(n,\mathtt{i},\mathtt{k},\mathtt{n}_i,\gamma_i,h)$
shown on the diagonal of the following commutative diagram
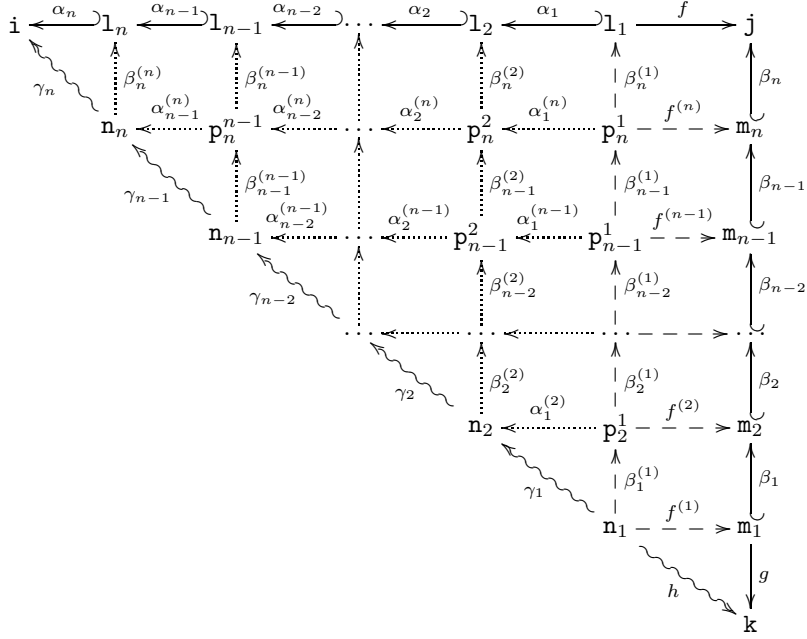
\begin{equation}\label{nice}
\xymatrix{
\mathtt{i} &
\mathtt{l}_n\ar@{_{(}->}[l]_{\alpha_n}
&\mathtt{l}_{n-1}\ar@{_{(}->}[l]_{\alpha_{n-1}}
&\dots\ar@{_{(}->}[l]_{\alpha_{n-2}}
&\mathtt{l}_2\ar@{_{(}->}[l]_{\alpha_{2}}
&\mathtt{l}_1\ar@{_{(}->}[l]_{\alpha_{1}}\ar[r]^f
&\mathtt{j}\\
&\mathtt{n}_n\ar@{.>}[u]_{\beta_n^{(n)}}
\ar@{~>}[ul]^{\gamma_n}
&\mathtt{p}_n^{n-1}\ar@{.>}[u]_{\beta_n^{(n-1)}}
\ar@{.>}[l]_{\alpha_{n-1}^{(n)}}
&\dots\ar@{.>}[l]_{\alpha_{n-2}^{(n)}}\ar@{.>}[u]
&\mathtt{p}_n^{2}\ar@{.>}[u]_{\beta_n^{(2)}}
\ar@{.>}[l]_{\alpha_{2}^{(n)}}
&\mathtt{p}_n^{1}\ar@{-->}[u]_{\beta_n^{(1)}}
\ar@{.>}[l]_{\alpha_{1}^{(n)}}\ar@{-->}[r]^{f^{(n)}}&
\mathtt{m}_n\ar@{_{(}->}[u]_{\beta_n}\\
&&\mathtt{n}_{n-1}\ar@{.>}[u]_{\beta_{n-1}^{(n-1)}}
\ar@{~>}[ul]^{\gamma_{n-1}}
&\dots\ar@{.>}[l]_{\alpha_{n-2}^{(n-1)}}\ar@{.>}[u]
&\mathtt{p}_{n-1}^{2}\ar@{.>}[l]_{\alpha_{2}^{(n-1)}}
\ar@{.>}[u]_{\beta_{n-1}^{(2)}}
&\mathtt{p}_{n-1}^{1}\ar@{.>}[l]_{\alpha_{1}^{(n-1)}}
\ar@{-->}[u]_{\beta_{n-1}^{(1)}}\ar@{-->}[r]^{f^{(n-1)}}&
\mathtt{m}_{n-1}\ar@{_{(}->}[u]_{\beta_{n-1}}\\
&&&\dots\ar@{~>}[ul]^{\gamma_{n-2}}\ar@{.>}[u]&
\dots\ar@{.>}[u]_{\beta_{n-2}^{(2)}}\ar@{.>}[l]
&\dots\ar@{-->}[u]_{\beta_{n-2}^{(1)}}\ar@{.>}[l]\ar@{-->}[r]&
\dots\ar@{_{(}->}[u]_{\beta_{n-2}}&\\
&&&&\mathtt{n}_2\ar@{~>}[ul]^{\gamma_{2}}
\ar@{.>}[u]_{\beta_{2}^{(2)}}
&\mathtt{p}_2^1\ar@{-->}[u]_{\beta_{2}^{(1)}}
\ar@{.>}[l]_{\alpha_{1}^{(2)}}\ar@{-->}[r]^{f^{(2)}}
&\mathtt{m}_2\ar@{_{(}->}[u]_{\beta_{2}}\\
&&&&&\mathtt{n}_1\ar@{~>}[ul]^{\gamma_{1}}
\ar@{-->}[u]_{\beta_{1}^{(1)}}\ar@{-->}[r]^{f^{(1)}}
\ar@{~>}[rd]_h
&\mathtt{m}_1\ar@{_{(}->}[u]_{\beta_{1}}\ar[d]^g\\
&&&&&&\mathtt{k}
}
\end{equation}
where the right column consists of pullbacks and all other
small squares are complement diagrams.
\end{proposition}

\begin{proof}
Follows from the construction of pullbacks in the category 
$\mathscr{P}(\mathscr{C})$ (see the proof of Theorem~\ref{thm81})
by induction on $n$.
\end{proof}

\subsection{Natural transformations}\label{s6.3}

Let $\mathscr{C}\in\mathscr{S}$ and $\mathtt{i},\mathtt{j}\in
\mathscr{C}$. For $n\geq 1$ and $s\in\{1,2,\dots,n+1\}$ we define 
the inclusion 
$\mathfrak{f}_n^{(s)}:\mathscr{P}^{n}(\mathtt{i},\mathtt{j})
\to\mathscr{P}^{n+1}(\mathtt{i},\mathtt{j})$ by 
mapping the element 
$\mathrm{D}(n,\mathtt{i},\mathtt{j},\mathtt{k}_i,\alpha_i,f)$ to
the element
\begin{displaymath}
\xymatrix{
&&&&&&&&\mathtt{j}\\
\mathtt{i}&\mathtt{k}_n\ar@{_{(}->}[l]_{\alpha_{n}}&
\dots\ar@{_{(}->}[l]_{\alpha_{n-1}}
&\mathtt{k}_{s+1}\ar@{_{(}->}[l]_{\alpha_{s+1}}&
\mathtt{k}_{s}\ar@{_{(}->}[l]_{\alpha_{s}}&
\mathtt{k}_{s}\ar@{=}[l]&
\mathtt{k}_{s-1}\ar@{_{(}->}[l]_{\alpha_{s-1}}&
\dots \ar@{_{(}->}[l]_{\alpha_{s-2}}&
\mathtt{k}_{1}\ar[u]^f\ar@{_{(}->}[l]_{\alpha_1}
},
\end{displaymath}
where for convenience we put  $\mathtt{k}_{n+1}:=\mathtt{i}$.

Further, for each $n\geq 2$ and each $s$, $2\leq s\leq n$, we 
define the surjection
$\mathfrak{t}_n^{(s)}:\mathscr{P}^{(n)}(\mathtt{i},\mathtt{j})
\to\mathscr{P}^{(n-1)}(\mathtt{i},\mathtt{j})$ by 
mapping the element 
$\mathrm{D}(n,\mathtt{i},\mathtt{j},\mathtt{k}_i,\alpha_i,f)$ to
the element
\begin{displaymath}
\xymatrix{
&&&&&&&\mathtt{j}\\
\mathtt{i}&\mathtt{k}_n\ar@{_{(}->}[l]_{\alpha_{n}}&
\dots\ar@{_{(}->}[l]_{\alpha_{n-1}}
&\mathtt{k}_{s+1}\ar@{_{(}->}[l]_{\alpha_{s+1}}&&
\mathtt{k}_{s-1}\ar@{_{(}->}[ll]_{\alpha_{s}\alpha_{s-1}}&
\dots \ar@{_{(}->}[l]_{\alpha_{s-2}}&
\mathtt{k}_{1}\ar[u]^f\ar@{_{(}->}[l]_{\alpha_1}
}
\end{displaymath}

We have the following statement:

\begin{proposition}\label{prom95}
Let $s$ and $n$ be as above.
\begin{enumerate}[(i)]
\item\label{prom95.1} $\mathfrak{f}_n^{(s)}$ is an injective
natural transformation from $\mathscr{P}^{n}$ to
$\mathscr{P}^{n+1}$.
\item\label{prom95.2} $\mathfrak{t}_n^{(s)}$ is an surjective
natural transformation from $\mathscr{P}^{n}$ to
$\mathscr{P}^{n-1}$.
\item\label{prom95.3} Both 
$\mathfrak{t}_{n+1}^{(s)}\mathfrak{f}_n^{(s)}$ and
$\mathfrak{t}_{n+1}^{(s+1)}\mathfrak{f}_n^{(s)}$ are
the identity transformations.
\item\label{prom95.4}  
$\mathfrak{t}_{n+1}^{(r)}\mathfrak{f}_n^{(s)}=
\mathfrak{f}_{n-1}^{(s)}\mathfrak{t}_{n}^{(r-1)}$ if 
$n\geq 2$ and $s+1<r$.
\item\label{prom95.5}  
$\mathfrak{t}_{n+1}^{(r)}\mathfrak{f}_n^{(s)}=
\mathfrak{f}_{n-1}^{(s-1)}\mathfrak{t}_{n}^{(r)}$ if 
$n\geq 2$ and $r<s$.
\end{enumerate}
\end{proposition}

We remark that for the quasi-iteration $\mathscr{P}^{(n)}$ there
is no analogue of $\mathfrak{f}_{n}^{(n+1)}$.

\begin{proof}
To prove \eqref{prom95.1} we have to show that 
$\mathfrak{f}_n^{(s)}$ behaves well with respect to 
the composition of morphisms. The latter one is described by 
the diagram \eqref{nice} (Proposition~\ref{prop91}).
Applying $\mathfrak{f}_n^{(s)}$ we double one row and one
column in \eqref{nice} (inserting the equality signes between
the doubled elements). The claim \eqref{prom95.1} would follow
if we would show that the obtained diagram is again of the
form \eqref{nice}. This reduces \eqref{prom95.1}
to the following facts:

\begin{lemma}\label{l92-1}
\begin{enumerate}[(a)]
\item\label{l92-1.1}
For any $\mathtt{i},\mathtt{j}\in\mathscr{C}$ and
$f\in \mathscr{C}(\mathtt{i},\mathtt{j})$ the diagrams
\begin{displaymath}
\xymatrix{
\mathtt{i}\ar[r]^f\ar@{=}[d]&\mathtt{j}\ar@{=}[d]\\
\mathtt{i}\ar[r]^f&\mathtt{j}\\
}\quad\text{and}\quad
\xymatrix{
\mathtt{j}&\mathtt{j}\ar@{=}[l]\\
\mathtt{i}\ar[u]^f&\mathtt{i}\ar[u]^f\ar@{=}[l]\\
}
\end{displaymath}
are pullbacks.
\item\label{l92-1.2}
For any $\mathtt{i},\mathtt{j}\in\mathscr{C}$ and any monomorphsism
$\alpha\in \mathscr{C}(\mathtt{i},\mathtt{j})$ the diagrams
\begin{displaymath}
\xymatrix{
\mathtt{i}\ar@{^{(}->}[r]^{\alpha}\ar@{=}[d]&\mathtt{j}\ar@{=}[d]\\
\mathtt{i}\ar@{^{(}->}[r]^{\alpha}&\mathtt{j}\\
}\quad\text{ and }\quad
\xymatrix{
\mathtt{j}&\mathtt{j}\ar@{=}[l]\\
\mathtt{i}\ar@{^{(}->}[u]^{\alpha}&
\mathtt{i}\ar@{^{(}->}[u]^{\alpha}\ar@{=}[l]
}
\end{displaymath}
are complement diagrams to the diagrams 
\begin{displaymath}
\xymatrix{\mathtt{i}\ar@{^{(}->}[r]^{\alpha}&\mathtt{j}\ar@{=}[r]
&\mathtt{j}}\quad\text{and}\quad
\xymatrix{\mathtt{i}\ar@{=}[r]&\mathtt{i}\ar@{^{(}->}[r]^{\alpha}
&\mathtt{j}}
\end{displaymath}
respectively.
\end{enumerate}
\end{lemma}

\begin{proof}
The whole statement \eqref{l92-1.1} and the statement \eqref{l92-1.2} for
the right diagram are obvious. To prove \eqref{l92-1.2} for the
left digram consider the solid part of the commutative diagram
\begin{displaymath}
\xymatrix{
\mathtt{k}\ar@{^{(}->}[rr]^{\beta} \ar@{-->}[rd]^{\gamma\delta^{-1}}
&& \mathtt{j}\ar@{=}[dd]\\
&\mathtt{i}\ar@{^{(}->}[ru]^{\alpha}\ar@{=}[d]&\\
\mathtt{l}\ar@{^{(}->}[r]^{\gamma}\ar@{^{(}->}[uu]^{\delta}
&\mathtt{i}\ar@{^{(}->}[r]^{\alpha}&\mathtt{j}\\
}
\end{displaymath}
and assume that the outer square is a pullback. Because of 
the latter assumption and \eqref{l92-1.1} we get that 
$\delta$ is an isomorphism and hence we have the induced
map $\gamma\delta^{-1}$ as required by Condition~\ref{pbm2}.
This completes the proof.
\end{proof}

The statement \eqref{prom95.2} is more complicated. We again
have to show that  $\mathfrak{t}_k^{(s)}$ behaves well with 
respect to the composition of morphisms.  Applying 
$\mathfrak{t}_k^{(s)}$ to \eqref{nice} just forgets one
row and one column in the diagram \eqref{nice}. The claim 
\eqref{prom95.2} would follow if we would show that the 
obtained diagram is again of the form \eqref{nice}. This 
reduces \eqref{prom95.2} to \cite[Proposition~7.2]{Mi} and the 
following statement:

\begin{lemma}\label{llA}
Assume that all small squares on the following diagrams are
complement diagrams:
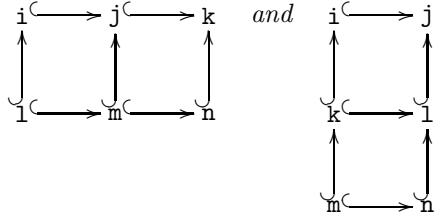
\begin{displaymath}
\xymatrix{
\mathtt{i}\ar@{^{(}->}[r]&\mathtt{j}\ar@{^{(}->}[r]&\mathtt{k}\\
\mathtt{l}\ar@{^{(}->}[r]\ar@{^{(}->}[u]
&\mathtt{m}\ar@{^{(}->}[r]\ar@{^{(}->}[u]
&\mathtt{n}\ar@{^{(}->}[u]
}\quad\text{and}\quad
\xymatrix{
\mathtt{i}\ar@{^{(}->}[r]&\mathtt{j}\\
\mathtt{k}\ar@{^{(}->}[r]\ar@{^{(}->}[u]&\mathtt{l}\ar@{^{(}->}[u] \\
\mathtt{m}\ar@{^{(}->}[r]\ar@{^{(}->}[u]&\mathtt{n}\ar@{^{(}->}[u] 
}
\end{displaymath}
Then the outer rectangles of both diagrams are complement
diagrams as well.
\end{lemma}

\begin{proof}
We start from the left diagram (which is easier to deal with). 
Consider  the extended diagram
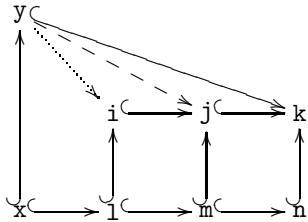
\begin{displaymath}
\xymatrix{
\mathtt{y}\ar@{^{(}->}[rrrd]\ar@{-->}[rrd]\ar@{.>}[rd]&&&\\
&\mathtt{i}\ar@{^{(}->}[r]&\mathtt{j}\ar@{^{(}->}[r]&\mathtt{k}\\
\mathtt{x}\ar@{^{(}->}[r]\ar@{^{(}->}[uu]
&\mathtt{l}\ar@{^{(}->}[r]\ar@{^{(}->}[u]
&\mathtt{m}\ar@{^{(}->}[r]\ar@{^{(}->}[u]
&\mathtt{n}\ar@{^{(}->}[u]
}
\end{displaymath}
such that the solid part commutes and the outer square is a
pullback. Since the right small square is a complement
diagram, we obtain the dashed map, which is automatically a
monomorphism. The condition that the outer square is a
pullback implies that the square with $\mathtt{x}$,
$\mathtt{y}$, $\mathtt{j}$ and $\mathtt{m}$ is a pullback as well.
Now since the middle small square is a complement
diagram, we obtain the dotted map, making the whole diagram commutative.
This implies that the outer rectangle of our two small squares
is a complement diagram.

Now we go on to the right diagram. Consider the following
diagram:
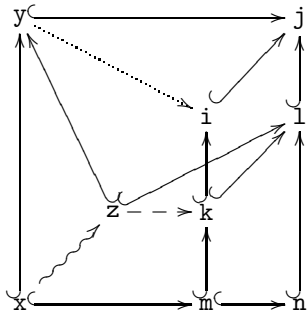
\begin{displaymath}
\xymatrix{
\mathtt{y}\ar@{^{(}->}[rrr]\ar@{.>}[rrd] &&&\mathtt{j}\\
&&\mathtt{i}\ar@{^{(}->}[ru]&\mathtt{l}\ar@{^{(}->}[u]\\
&\mathtt{z}\ar@{^{(}->}[rru]\ar@{_{(}->}[luu]\ar@{-->}[r]&
\mathtt{k}\ar@{^{(}->}[ru]\ar@{^{(}->}[u]&\\
\mathtt{x}\ar@{^{(}->}[rr]\ar@{~>}[ru]\ar@{^{(}->}[uuu]
&&\mathtt{m}\ar@{^{(}->}[r]\ar@{^{(}->}[u]&
\mathtt{n}\ar@{^{(}->}[uu]
}
\end{displaymath}
Assume that the outer square is a pullback and that the
square containing $\mathtt{z}$, $\mathtt{y}$, $\mathtt{l}$
and $\mathtt{j}$ is a pullback as well. The latter condition
implies the existence of the curled map from $\mathtt{x}$
to $\mathtt{z}$ as indicated. Since the outer square is a pullback,
it follows that the square containing $\mathtt{x}$, $\mathtt{z}$, 
$\mathtt{l}$ and $\mathtt{n}$ is a pullback as well.
Hence we can use the fact that the right bottom small
square is a complement diagram and obtain a dashed map 
(in fact a monomorphism) from
$\mathtt{z}$ to $\mathtt{k}$ as indicated and everything commutes.
It is left to observe that the right top small square is a
complement diagram, and hence there should exist the dotted map
from $\mathtt{y}$ to $\mathtt{i}$ making the whole diagram 
commutative. The necessary claim follows.
\end{proof}

The statements \eqref{prom95.3}--\eqref{prom95.5} are proved by
a direct calculation.
\end{proof}

Altogether we have the following picture:

\begin{equation}\label{eq91-1}
\xymatrix{
\mathscr{P}^{0}\ar@/_/@{^{(}->}[r]|-{\mathfrak{f}_0^1} &
\mathscr{P}^{1}\ar@/_/@{^{(}->}[r]|-{\mathfrak{f}_1^1}
\ar@/_1.5pc/@{^{(}->}[r]|-{\mathfrak{f}_1^2} &
\mathscr{P}^{2}\ar@/_/@{^{(}->}[r]|-{\mathfrak{f}_2^1}
\ar@/_1.5pc/@{^{(}->}[r]|-{\mathfrak{f}_2^2}
\ar@/_2.5pc/@{^{(}->}[r]|-{\mathfrak{f}_2^3}
\ar@/_/@{->>}[l]|-{\mathfrak{t}_2^2}
&\mathscr{P}^{3}\ar@/_/@{^{(}->}[r]|-{\mathfrak{f}_3^1}
\ar@/_1.5pc/@{^{(}->}[r]|-{\mathfrak{f}_3^2}
\ar@/_2.5pc/@{^{(}->}[r]|-{\mathfrak{f}_3^3}
\ar@/_3.5pc/@{^{(}->}[r]|-{\mathfrak{f}_3^4}
\ar@/_/@{->>}[l]|-{\mathfrak{t}_3^2}
\ar@/_1.5pc/@{->>}[l]|-{\mathfrak{t}_3^3}&\dots
\ar@/_/@{->>}[l]|-{\mathfrak{t}_4^2}
\ar@/_1.5pc/@{->>}[l]|-{\mathfrak{t}_4^3}
\ar@/_2.5pc/@{->>}[l]|-{\mathfrak{t}_4^4}
}
\end{equation}

Let $\mathscr{A}$ denote the category, whose objects are
$\mathscr{P}^{n}$, $n=0,1,2\dots$, and morphisms are all
possible natural transformations of functors.

Let $\mathscr{B}$ denote the category, whose objects are
$n\in\mathbb{N}$, and for $m,n\in \mathscr{B}$ we have
\begin{displaymath}
\mathscr{B}(m,n)=
\left\{
\left(
\begin{array}{ccccc}
1&2&\dots&n&\infty\\
a_1&a_2&\dots&a_n&a_{\infty}
\end{array}
\right):
\begin{array}{ll}
a_1=1;& a_i\leq a_{i+1};\\
a_{\infty}=\infty;& a_i\in\{1,\dots,m,\infty\}
\end{array}
\right\}
\end{displaymath}
with the obvious multiplication. One shows that any morphisms
in the category $\mathscr{B}$ can be written as a composition
of the following morphisms:
\begin{displaymath}
\varphi_n^{(s)}=
\left(
\begin{array}{ccccccccccc}
1&2&\dots&s-1&s&s+1&s+2&\dots&n&n+1&\infty\\
1&2&\dots&s-1&s&s&s+1&\dots&n-1&n&\infty\\
\end{array}
\right),
\end{displaymath}
where $n\in\mathbb{N}$ and $s\in\{1,2,\dots,n\}$; 
\begin{displaymath}
\varphi_n^{(n+1)}=
\left(
\begin{array}{cccccc}
1&2&\dots&n&n+1&\infty\\
1&2&\dots&n&\infty&\infty\\
\end{array}
\right)
\end{displaymath}
and
\begin{displaymath}
\tau_n^{(s)}=
\left(
\begin{array}{cccccccccc}
1&2&\dots&s-1&s&s+1&\dots&n-1&n&\infty\\
1&2&\dots&s-1&s+1&s+2&\dots&n&n+1&\infty\\
\end{array}
\right),
\end{displaymath}
where $n>1$ and $s\in\{2,\dots,n\}$.
Denote by $\mathcal{O}_{n}$ the monoid of all 
order-preserving transformations on the chain 
$\{1,\dots,n,\infty\}$, see for example \cite{Gl}, and
by $\mathcal{O}'_{n}$ the submonoid of 
$\mathcal{O}_{n}$ consisting of all all 
transformations, which fix the points $1$ and $\infty$. 
One easily shows that for $n\in \mathbb{N}$ the monoid 
$\mathrm{End}_{\mathscr{B}}(n)$ is isomorphic to 
$\mathcal{O}'_{n}$.

\begin{proposition}\label{pr95-1}
The assignment $n\mapsto \mathscr{P}^n$,
$\varphi_n^{(s)}\mapsto \mathfrak{f}_n^{(s)}$, and
$\tau_n^{(s)}\mapsto \mathfrak{t}_n^{(s)}$, extends to a 
faithful functor, $\mathrm{F}:\mathscr{B}\to \mathscr{A}$.
\end{proposition}

\begin{proof}
A direct calculation shows that $\mathrm{F}$ exists.
Let $\varphi,\psi\in \mathscr{B}(m,n)$. To claim that $F$ 
is  faithfull it is enough to find $\mathscr{C}\in 
\mathscr{S}$ and  $\mathtt{i}\in \mathscr{C}$ such 
that  $F(\varphi)$ and $F(\psi)$ induce different 
morphisms from $\mathrm{End}_{\mathscr{P}^m(\mathscr{C})}(\mathtt{i})$
to $\mathrm{End}_{\mathscr{P}^n(\mathscr{C})}(\mathtt{i})$. 
Later on in Subsection~\ref{s6.4} we will show that 
the category $\mathscr{C}_1$ from Subsection~\ref{s3.1}
belongs to $\mathscr{S}$. The property above is then
easily verified by a direct calculation if one
takes $\mathtt{i}$ to be a finite set of cardinality at
least $\max(m,n)+1$.
\end{proof}

As an immediate corollary we obtain:

\begin{corollary}\label{cor9605}
Let $\mathscr{C}\in \mathscr{S}$ and  $\mathtt{i}\in \mathscr{C}$. 
Then the functor $F$ induces an action of the monoid
$\mathcal{O}'_{n-1}$ on the monoid
$\mathrm{End}_{\mathscr{P}^n(\mathscr{C})}(\mathtt{i})$  
by endomorphisms.
\end{corollary}

\subsection{Connection between $\mathscr{P}^{(n)}$ and
$\mathscr{P}^{n}$}\label{s6.35}

As we have seen, for $\mathscr{C}\in\mathscr{S}$ the morphism
sets in $\mathscr{P}^{(n)}(\mathscr{C})$ and
$\mathscr{P}^{n}(\mathscr{C})$ can be canonically identified.
However, the products are rather different. Nevertheless,
there is a clear connection between them, which can be 
described as follows: Let us for the moment denote the product 
in $\mathscr{P}^{(n)}(\mathscr{C})$ by $*$ and the product 
in $\mathscr{P}^{n}(\mathscr{C})$ by $\star$. We have the 
following:

\begin{proposition}\label{prop97}
Let $x=\mathrm{D}(n,\mathtt{i},\mathtt{j},\mathtt{l}_i,\alpha_i,f)$ and 
$y=\mathrm{D}(n,\mathtt{j},\mathtt{k},\mathtt{m}_i,\beta_i,g)$. 
Set $\mathfrak{a}=
\mathfrak{f}_{n-1}^{(1)}\circ\dots\circ\mathfrak{f}_2^{(1)}\circ
\mathfrak{f}_1^{(1)}\circ\mathfrak{t}_2^{(2)}\circ\dots\circ
\mathfrak{t}_{n-1}^{(n-1)}\circ\mathfrak{t}_n^{(n)}$.
Then $x*y=x\star \mathfrak{a}(y)$.
\end{proposition}

\begin{proof}
Follows from Proposition~\ref{prop91}
and Lemma~\ref{l92-1} by a direct calculation.
\end{proof}

\begin{remark}
{\rm The statement of Proposition~\ref{prop97} reminds of the
following fact: if $(S,\star)$ is a semigroup and $\sigma$ is an 
idempotent  endomorphism of $S$ (retraction), then 
$(S,*)$, where  $x*y:=x\star \sigma(y)$, is a semigroup.
}
\end{remark}

\subsection{$\mathscr{P}^n$ and orthodox semigroups}\label{s6.38}

\begin{theorem}\label{thm992}
Let $\mathscr{C}\in \mathscr{S}$ be such that 
all morphisms in $\mathscr{C}$ are monomorphisms.
\begin{enumerate}[(i)]
\item \label{thm992.01}
Let $x=\mathrm{D}(n,\mathtt{i},\mathrm{j},\mathrm{k}_i,\alpha_i,f)$
and $y=\mathrm{D}(n,\mathtt{j},\mathrm{i},\mathrm{l}_i,\beta_i,g)$.
Then $x$ and $y$ form a pair of inverse elements if and only if
there is an isomorphism $\gamma:\mathrm{k}_1\to \mathrm{l}_1$ such that
$\alpha_n\dots \alpha_2\alpha_1=g\gamma$ and
$f=\beta_n\dots \beta_2\beta_1\gamma$.
\item \label{thm992.1}
For each $k\geq 0$ and for each $\mathtt{i}\in \mathscr{C}$ the
monoid $\mathrm{End}_{\mathscr{P}^{k}(\mathscr{C})}
\big(\mathtt{i}\big)$ is regular and is a retract of the 
monoid $\mathrm{End}_{\mathscr{P}^{k+1}(\mathscr{C})}
\big(\mathtt{i}\big)$.
\item \label{thm992.2}
The element $\mathrm{D}(n,\mathtt{i},\mathrm{j},\mathrm{k}_i,\alpha_i,f)$
is an idempotent if and only if $f=\alpha_n\dots\alpha_2\alpha_1$.
\item \label{thm992.3}
For each $k\geq 0$ and for each $\mathtt{i}\in \mathscr{C}$ the
monoid $\mathrm{End}_{\mathscr{P}^{k}(\mathscr{C})}
\big(\mathtt{i}\big)$ is orthodox
(i.e. it is regular and its idempotents form a subsemigroup).
\end{enumerate}
\end{theorem}

\begin{proof}
One shows by a direct calculation (using Proposition~\ref{prop91} 
and Lemma~\ref{l92-1}) that any pair of elements satisfying the
conditions of \eqref{thm992.01} is inverse to each other.
Let $x$ and $y$ be a pair of inverse elements. Consider the 
retraction $\mathfrak{a}$  from Proposition~\ref{prop97}.  
Obviously, $\mathfrak{a}(x)$ and $\mathfrak{a}(y)$ constitute a 
pair of inverse elements as well. Now the claim \eqref{thm992.01} 
follows from Lemma~\ref{lem:regularity} and 
Proposition~\ref{pmonoinv} (note that in an inverse semigroup 
the inverse element is unique).

The statement \eqref{thm992.1} now follows from
\eqref{thm992.01} and Proposition~\ref{prom95}\eqref{prom95.3}.

That all elements having the form as in \eqref{thm992.2}
are idempotents is checked by a direct calculation. Let 
$x$ be an idempotent. Consider the retraction $\mathfrak{a}$ 
from Proposition~\ref{prop97}.  Obviously, $\mathfrak{a}(x)$ 
is an idempotent as well. Now the claim \eqref{thm992.2} 
follows from Lemma~\ref{lem:idempotents}.

From \eqref{thm992.01} and \eqref{thm992.2} it follows that
each inverse of an idempotent in the semigroup 
$\mathrm{End}_{\mathscr{P}^{k}(\mathscr{C})}
\big(\mathtt{i}\big)$ is an idempotent itself.
Hence \eqref{thm992.3} follows from \eqref{thm992.1} and
\cite[Theorem~1.1]{How}.
\end{proof}

\subsection{Applications to finite sets}\label{s6.4}

\begin{proposition}\label{prop991}
$\mathscr{C}_1\in\mathscr{S}$ (see Subsection~\ref{s3.1}).
\end{proposition}

\begin{proof}
From Subsection~\ref{s3.1} we know that $\mathscr{C}_1$ satisfies
Condition~\ref{pbm}. Hence we have only to check that
$\mathscr{C}_1$ satisfies Condition~\ref{pbm2}.
Let $X,Y,Z$ be finite sets, $\alpha:X\hookrightarrow Y$
and $\beta:Y\hookrightarrow Z$. Set $U=\beta(\alpha(X)) \cup 
(Z\setminus \beta(Y))$. We have the natural inclusion
$\mathrm{incl}:U\hookrightarrow Z$. We also define
the map $\gamma:X\to U$ via $\gamma(x)=\beta(\alpha(x))$.
Then $\gamma$ in obviously an inclusion. A direct calculation
shows that the diagram
\begin{displaymath}
\xymatrix{
U\ar@{^{(}->}[r]^{\mathrm{incl}}&Z\\
X\ar@{^{(}->}[r]^{\alpha}\ar@{^{(}->}[u]^{\gamma}
&Y\ar@{^{(}->}[u]^{\beta}
}
\end{displaymath}
is a complement diagram. The claim follows.
\end{proof}

\begin{remark}
{\rm 
The category $\mathscr{C}_2^{\mathrm{op}}$ from 
Subsection~\ref{s3.15} satisfies Condition~\ref{pbm}. One
can show that it does not satisfy Condition~\ref{pbm2}.
}
\end{remark}

Proposition~\ref{prop991} allows us to consider the categories
$\mathscr{P}^k(\mathscr{C}_1)$ for each $k\in\mathbb{N}$. In 
particular, for each $n\in \mathbb{N}$ we have an orthodox semigroup
\begin{displaymath}
\mathcal{RS}(n,k):=
\mathrm{End}_{\mathscr{P}^k(\mathscr{C}_1)}(\mathtt{n}),
\end{displaymath}
where $\mathtt{n}=\{1,2,\dots,n\}$. We have
$\mathcal{RS}(n,1)\cong\mathcal{IS}_n$, the classical
symmetric inverse semigroup. In what follows we list some basis 
properties of $\mathcal{RS}(n,k)$ (the proofs are left to the reader).

\begin{enumerate}[(I)]
\item\label{thmrslast.1}
Each element of $\mathcal{RS}(n,k)$ can be uniquely written in
the form
\begin{equation}
\mathrm{D}(k,A_i,f):\quad
\xymatrix{
\mathtt{n} & A_k\ar@{_{(}->}[l]_{\mathrm{incl}}
&\dots\ar@{_{(}->}[l]_{\mathrm{incl}}
 & A_2\ar@{_{(}->}[l]_{\mathrm{incl}} & 
A_1\ar@{_{(}->}[l]_{\mathrm{incl}}\ar[r]^f & \mathtt{n} 
},
\end{equation}
where $A_1\subset A_2\subset \dots\subset A_k\subset \mathtt{n}$,
$\mathrm{incl}$ denotes the natural inclusion, and 
$f:A_1\to \mathtt{n}$ is an injection.
\item\label{thmrslast.2}
$|\mathcal{RS}(n,k)|=\sum_{i=0}^n\binom{n}{i}^2i!k^{n-i}$.
\item\label{thmrslast.3}
$|E(\mathcal{RS}(n,k))|=(k+1)^{n}$.
\item\label{thmrslast.4}
$\mathrm{D}(k,A_i,f)\,\mathcal{D}\,\mathrm{D}(k,B_i,g)$ if
and only if $\mathrm{D}(k,A_i,f)\,\mathcal{J}\,\mathrm{D}(k,B_i,g)$ 
if and only if $|A_1|=|B_1|$.
\item\label{thmrslast.5}
$\mathrm{D}(k,A_i,f)\,\mathcal{R}\,\mathrm{D}(k,B_i,g)$ if
and only if $\mathrm{Im}(f)=\mathrm{Im}(g)$. The $\mathcal{R}$-class 
of the element $\mathrm{D}(k,A_i,f)$ contains 
$k^{n-|A_1|}$ idempotents.
\item\label{thmrslast.6}
$\mathrm{D}(k,A_i,f)\,\mathcal{L}\,\mathrm{D}(k,B_i,g)$ if
and only if $(A_i)=(B_i)$. Each $\mathcal{L}$-class of 
$\mathcal{RS}(n,k)$ contains a unique idempotent.
\item\label{thmrslast.7}
$\mathrm{D}(k,A_i,f)\,\mathcal{H}\,\mathrm{D}(k,B_i,g)$ 
if and only if $\mathrm{Im}(f)=\mathrm{Im}(g)$
and $(A_i)=(B_i)$.
\item\label{thmrslast.8}
All maximal subgroups of $\mathcal{RS}(n,k)$ are isomorphic
to symmetric groups of rank $\leq n$.
\item\label{thmrslast.9}
$E(\mathcal{RS}(n,k))$ is a Boolean of right singular semigroups,
i.e. there is an epimorphism (induced by $\mathfrak{a}$
from Proposition~\ref{prop97}), whose 
image is the Boolean $(2^{\mathtt{n}},\cap)$ and such that
each congruence class of the kernel is a semigroup of right zeros.
The kernel of this epimorphism coincides with the minimum
semilatice congruence on $E(\mathcal{RS}(n,k))$. 
\item\label{thmrslast.10}
By \eqref{thmrslast.1} and Theorem~\ref{thm992}\eqref{thm992.2}
idempotents in $\mathcal{RS}(n,k)$ are described by flags
$(A_1,\dots,A_k)$ of subsets of $\mathtt{n}$. In this notation, the
multiplication of idempotents is as follows:
\begin{displaymath}
(B_1,\dots,B_k)\cdot(A_1,\dots,A_k)=
\big((A_1\cap B_i)\cup(A_i\setminus A_1)\big)_{i=1}^k.
\end{displaymath}
\end{enumerate}

\begin{remark}\label{lastrem}
{\rm 
Consider the category $\mathscr{C}_{10}$ from 
Subsection~\ref{s3.15new}. One shows that $\mathscr{C}_{10}$
satisfies Condition~\ref{pbm2}. One further easily computes that
$\mathrm{End}_{\mathscr{P}^k(\mathscr{C}_{10})}(\mathbb{N})$
is a bisimple orthodox monoid for each $k\geq 1$ (it is inverse
for $k=1$).
}
\end{remark}

\end{document}